\documentclass[11pt,reqno]{amsart}

\usepackage{graphicx}
\usepackage{amsmath}
\usepackage{amssymb}
\usepackage{amsfonts}
\usepackage{amsthm}
\usepackage{epstopdf}

\usepackage{color}
\usepackage{eucal}

\usepackage{caption}
\usepackage{subcaption}

\usepackage{wrapfig}
\usepackage{floatflt}

\usepackage{comment}

\newcommand{\R}{\mathbb R}
\newcommand{\dd}{\,d}

\newcommand{\ds}{\displaystyle}

\newtheorem{conjecture}{Conjecture} 
\theoremstyle{remark}

\newcommand{\defeq}{\stackrel{\rm{def}}{=}}

\begin{document}
\title[3D critical Zakharov-Kuznetsov equation]
{The 3D critical Zakharov--Kuznetsov equation: \\
blow-up and soliton dynamics}

\author[C. Klein]{Christian Klein}
\address{Université Bourgogne Europe, CNRS \\
Institut de Math\'ematiques de Bourgogne, UMR 5584 \\
 21000 Dijon, France;\\
 Institut Universitaire de France} 
\email{Christian.Klein@u-bourgogne.fr}

\author[S. Roudenko]{Svetlana Roudenko}
\address{Department of Mathematics \& Statistics\\Florida International University,  
Miami, FL, 33199, USA}
\curraddr{}
\email{sroudenko@fiu.edu}
    
\author[N. Stoilov]{Nikola Stoilov}
\address{Université Bourgogne Europe, CNRS \\
Institut de Math\'ematiques de Bourgogne, UMR 5584 \\
 21000 Dijon, France} 
\email{Nikola.Stoilov@u-bourgogne.fr}



\begin{abstract}
We study the full three-dimensional dynamics of the $L^2$-critical Zakharov-Kuznetsov equation with the fractional nonlinearity $|u|^{4/3}u$, equivalently $u^{7/3}$ for real-valued functions.  
This equation is a higher-dimensional extension of the generalized Korteweg-de Vries equation.  
In the critical setting solutions to this 3D ZK equation may blow up in finite time or exhibit global time dynamics. 
The novelties of this work is to treat a non-integer power and to study the dynamics of solutions in a higher dimension.   

We first review the finite time blow-up in 2D critical ZK, then do a formal analysis of the slightly mass-supercritical blow-up dynamics, deriving the corrections to the blow-up rate and profile for the critical ZK equation in any dimension.
We then perform a computational study of solutions, utilizing full 3D numerical simulations. 
In particular, we use a Fourier pseudospectral discretization and an 
integrating factor fourth-order Runge-Kutta method on a full three-dimensional grid. 
A multi-GPU implementation makes it possible to follow blow-up solutions in a full 3D setting. We examine perturbations of the ground state, Gaussian data, and nonsymmetric two-bump configurations. The computations show dispersive and concentrating regimes, in both cases with radiation emitted in a conic-type region opposite to the direction of propagation and convergence of the concentrating core toward a rescaled ground-state profile.  The two-bump experiments also demonstrate that total mass alone does not determine the blow-up dynamics.  We discuss the numerical evidence for the predicted blow-up rate and identify the pre-asymptotic and resolution limitations that remain near the blow-up time.
\end{abstract}

\maketitle

\section{Introduction}
In this paper we study 
the 3D critical Zakharov-Kuznetsov (ZK) equation
\begin{equation*}
u_t + (u_{xx}+u_{yy}+u_{zz} + |u|^{4/3}u)_x  = 0,
\end{equation*}
where $u=u(x,y,z,t)$ is real-valued, $(x,y,z) \in \mathbb R^3$, and $t \in \mathbb R$. Since we consider only real-valued solutions $u$, the absolute value in the nonlinearity can be dropped to simply write
\begin{equation}\label{ZK}
u_t + (u_{xx}+u_{yy}+u_{zz} + u^{7/3})_x  = 0.
\end{equation}
This equation is a three-dimensional variant of the generalized  
Korteweg-de Vries (gKdV) equation, and the square power $(u^2)_x$ is the celebrated (one-dimensional) KdV model for weakly nonlinear waves in shallow water. The square power is also the standard 3D ZK equation, which was originally proposed by Zakharov and Kuznetsov in the description of 
weakly magnetized ion-acoustic waves in a low-pressure magnetized 
plasma \cite{ZK1974}, where they raised the question of soliton stability  
in a higher-dimensional setting. While originally the equation was proposed by Zakharov and Kuznetsov in the 3D setting, the first rigorous derivation as a long-wave small-amplitude limit of the Euler-Poisson system in the cold-plasma approximation was done by Lannes, Linares and Saut in \cite{LLS}, see also \cite{LS}. Other derivations exist as well, for a review see \cite{LLS, FHRY2} and references therein.

During their lifespans, solutions $u(t)$ to the ZK equation conserve 
energy (Hamiltonian), $L^2$-norm (often called momentum or mass) and the integral, defined as follows
\begin{align}\label{MC}
M[u(t)] &\defeq\int_{\mathbb{R}^3} |u(t,x,y,z)|^2 \, dx dy dz = M[u(0)], \\
\label{EC}
E[u(t)] & \defeq \dfrac{1}{2}\int_{\mathbb{R}^3}
\big[u_{x}^{2}(t)+u_{y}^2(t)+u_z^2(t)\big] \, dx dy dz 
- \dfrac{3}{10}\int_{\mathbb{R}^3} |u(t)|^{\frac{10}3} \, dx dy dz 
= E[u(0)], \\
\label{L1-inv}
\int_{\mathbb{R}}  u(x, & \, y, z,t) \, dx = \int_{\mathbb{R}} u(x,y,z,0) \, dx.
\end{align}

The equation \eqref{ZK} has a scaling invariance: 
if $u(x,y,z,t)$ is a solution of \eqref{ZK}, then so is the rescaled version 
\begin{equation}\label{E:scaling}
u_\lambda(x,y,z,t) 
= \lambda^{3/2} u(\lambda x, \lambda y, \lambda z,\lambda^3 t), \quad \lambda > 0.
\end{equation}
This symmetry makes the Sobolev norm $\dot{H}^s$ with $s=0$ invariant, thus, making the equation \eqref{ZK} $L^2$-critical. 
The 3D critical ZK equation has other invariances such as translation in space and time.

The well-posedness for the Cauchy problem for the original 3D (quadratic) ZK equation with $H^s$ initial data has attracted significant interest in the last decade as it was a challenge to obtain it in the energy class: for a while the local well-posedness was only known in $H^{1+}(\mathbb R^3)$, see \cite{RV}, \cite{MP}. For solutions close to the solitary wave, it was shown in \cite{FHRY2} that a bootstrap argument can improve this and include $H^1$ solutions, and later, the work of Herr \& Kinoshita \cite{HK} established the local well-posedness in $H^s$ for $s>-\frac12$, and thus, sufficient for the $H^1$ energy solutions. Previous works include Linares \&  Saut \cite{LS} (local well-posedness in $H^s$ with $s>\frac98$), which was further improved by Ribaud \& Vento \cite{RV} (down to $H^s$ with $s>1$) and then extended to global well-posedness in  $H^s$, $s>1$ by Molinet \& Pilod \cite{MP}. 

For other nonlinearities in the 3D ZK equation, the well-posedness question is mostly open, except for the work of Linares \& Ramos \cite{LR2021}, where the authors establish well-posedness in $H^s$ for $s \in (\frac34,1)$ and an `almost well-posedness' for $s \in [1,2)$ (the data-to-solution map is continuous in a slightly larger space). 
In this work we are interested in $H^1$ solutions with an exponential 
decay (such as a sech or Gaussian), and thus, the above well-posedness is sufficient. It would be interesting to obtain a `full' well-posedness in this critical case to better understand this critical problem and investigate further the critical case in terms of the long term dynamics, as well as well-posedness for other nonlinearities of ZK equation in 3D.

In our paper \cite{KRS2} we investigated the 3D quadratic ZK equation, which is subcritical, and thus, all solutions are global,  we therefore confirmed numerically the soliton resolution conjecture. Furthermore, in 
\cite{FHRY2} the asymptotic stability of solitons was proved in the energy space $H^1$ for this equation, with further proof of the spectral properties for virial operator in \cite{HR2026}; see related works in 2D  \cite{KRS1},\cite{CMPS}, \cite{FHRY}, \cite{FHR3}, \cite{HR2025}. 
Thus, by now it is understood how solutions behave in the subcritical case (the original 3D quadratic ZK equation), however, 
as far as the $L^2$-critical ZK equation in 3D is concerned, there is not much work available, except for \cite{LR2021}. Meanwhile, it is important to investigate the $L^2$ critical case in 3D, especially, since it is a {\it non-integer} power in higher dimensions. This paper is a contribution to such investigations.  
\smallskip

The main goal of the present work is to investigate soliton formation, propagation, their instability leading to blow-up, radiation regime in the critical 3D ZK equation \eqref{ZK} from the numerical point of view and formal analysis. 
\smallskip

The equation has a family of localized traveling waves called solitary waves (often also referred to as solitons, although the model is not integrable), moving only in the positive $x$-direction:
\begin{equation}\label{E:wave}
u(x,y,z,t) = Q_c(x - c \, t, y, z), ~~c>0,
\end{equation}
where $Q_c$ is the dilation 
\begin{equation}\label{dilation}
Q_c(X) = c^{3/4} \, Q(\sqrt{c} \, X), \quad X = (x,y,z)
\end{equation}
We consider solitary waves that vanish at infinity, thus, $Q$ is the vanishing at infinity ground state solution of the nonlinear elliptic equation 
\begin{equation}\label{E:Q}
-Q + \Delta_{\mathbb R^3} Q   + Q^{\frac73} = 0,
\end{equation}  
i.e., the unique positive radial solution, up to translations, in (the energy space) $H^1(\mathbb{R}^3)$. This solution is strictly decreasing in the radial variable and exponentially decaying, see for example, \cite{K1989}. 
In particular, $Q \in C^{\infty}(\mathbb{R}^3)$, $\partial_r Q(r) <0$ for any $r = |(x,y,z)|>0$, and for any multi-index $\alpha$ there exists $\gamma>0$ such that
\begin{equation}\label{prop-Q}
|\partial^\alpha Q(X)| \lesssim_\alpha e^{-\gamma|X|} \quad \mbox{for any}\quad {X = (x,y,z)} \in \mathbb{R}^3.
\end{equation}

Since this is a critical case, it is expected that the ground state solutions are {\it not stable}, that is one of the questions we investigate in this paper. 
\bigskip

The first objective of this paper is to adapt formally to dimension three the next-order blow-up-profile mechanism obtained by Chen et al. in \cite{Gong} in the two-dimensional cubic ZK equation, which in its turn was inspired by results for the critical gKdV of Martel et al. in \cite{MMR}. This approach produces an explicit transverse-type quotient $\theta_{3d}$, see \eqref{E:theta_3}, and the predicted blow-up rate exponent
$$
\beta_{3d}=\frac1{3-\theta_{3d}}.
$$
We give two formulas suitable for accurate computation: the original transverse-type representation and a one-dimensional radial Fourier-slice formula, where we use the fact that the ground state is radial. To better motivate our analysis, we first provide a short review of the two-dimensional finite time blow-up in the 2D critical (modified) ZK equation, explaining the order of approximations that can be taken when studying the blow-up phenomena. We then formally derive the blow-up rates for {\it any} dimension, providing several ways, including radial, to compute the rate corrections. In particular, we compute the predicted blow-up rate for the 3D critical setting (as well as for the 2D and 4D).  
\smallskip

The second objective is computational.  We solve ~\eqref{ZK} on a full three-dimensional periodic box by a Fourier pseudospectral method and an integrating-factor RK4 scheme.  The $2^{12}\times2^9\times2^9$ discretization contains more than one billion spatial points, and the Fourier transforms are distributed over fully connected GPUs.  In our simulations, we investigate stability of solitons, soliton resolution, the radiation regimes, we track the formal profile and rate, compare with the cylindrically symmetric calculations in~\cite{KRS_cyl}, and explore nonsymmetric two-bump configurations that cannot be represented by a cylindrical reduction.

Investigating this non-integrable model, with non-integer power and in computationally-demanding higher dimension, our numerical approaches, in the light of the following conjecture, make definite progress towards providing numerical evidence for it:  

\begin{conjecture} \label{C:1} 
Solutions of the 3D critical ZK equation \eqref{ZK} with sufficiently smooth and single-peak localized initial data (e.g., $u_0 \in H^1$  and some sufficient decay) behave as follows:
\begin{enumerate}
\item[(a)] 
if $\|u_0\|_{L^2(\mathbb R^3)} < \|Q\|_{L^2(\mathbb R^3)}$, the solution $u(t)$ is global and  disperses into radiation, 

\item[(b)] 
if $\|u_0\|_{L^2(\mathbb R^3)} > \|Q\|_{L^2(\mathbb R^3)}$ and $E[u_0]<0$, then the solution $u(t)$ blows up in finite time\footnote{The finite time may depend on the decay at infinity of the initial data.}.
\end{enumerate}
In case (b) the solution decouples into a sum of rescaled concentrating solitons, moving in the positive $x$-direction, and radiation, dispersing away in a cone-type region opposite to the direction of propagation (negative $x$-direction); in case (a) only the radiative part is left. 
\end{conjecture}

Observe that in the $L^2$-critical case all solitary waves carry the same mass, since $\|Q_c\|_{L^2} = \|Q\|_{L^2}$ for every $c>0$ by \eqref{dilation}, and therefore, no soliton can be present below the threshold in part (a). We also point out that the assumption $E[u_0]<0$ in part (b) already implies $M[u_0]>M[Q]$, see the coercivity bound \eqref{E:coercive-lower} below, and hence, the mass condition in (b) is stated only for emphasis.

We note that the assumption about the {\it single-peak} localization is essential, since we show examples with two-bump initial data, having mass greater than that of the ground state, however, solutions have different outcomes, some radiate away and some blow up in finite time, see Sections \ref{S:doubleG} and \ref{S:doubleQ}.  

\medskip

The paper is organized as follows: in Section~\ref{S:rate} we review the two-dimensional rate law and give the formal derivation and computation of the three-dimensional analogue, as well as provide a formal rate correction for any dimension. In Section~\ref{S:ground} we record the ground-state identities and numerical diagnostics used for verification.  Section~\ref{S:Num} describes the full three-dimensional Fourier/GPU discretization. In Section~\ref{S:results} we present perturbations of the ground state and Gaussian data, followed by the nonsymmetric double-Gaussian and double-soliton experiments. We summarize the conclusions and the remaining numerical issues in Section~\ref{S:conclusions}.

\smallskip

{\bf Acknowledgements.} CK and NS were partially supported by the ANR project ISAAC-ANR-23-CE40-0015-01 and  the EIPHI Graduate School (contract ANR-17-EURE-0002). SR was partially supported by the NSF grant DMS-2452782.

\section{Formal blow-up rate and profile}\label{S:rate}
We first review the recent results on blow-up and its mechanism in the two-dimensional critical ZK equation and then transfer its leading-order algebra to the three-dimensional problem as well as to higher dimensions.  

\subsection{The two-dimensional critical ZK blow-up}\label{S:2D-rate}

In \cite{FHRY} the second author and collaborators proved that a solution with mass slightly supercritical and negative energy blows up in finite or infinite time, following the strategy of Merle \cite{M} for proving the existence of blow-up solutions in the critical gKdV equation. 
The proof in \cite{FHRY} did not assume any further decay on the initial data, which 
would be needed to conclude a finite time blow-up, since the decay on the far right tail is what gives the control of the finite time in the blow-up setting. 
They also introduced the crucial spectral virial step, which was later examined and proved with analytical methods and only verifying the signs of some inner products in \cite{HR2025}. 
The follow up work by Chen et al in \cite{Gong} used the approach of approximate solutions from \cite{MMR} as well as several key features from \cite{FHRY} such as monotonicity and the virial spectral property from \cite{FHRY, HR2025} to classify the behavior of solutions from a certain class of initial data with a prescribed decay, and in particular, proved the existence of the {\it finite time} blow-up. (This is similar to the 1D case of the critical gKdV equation, where in order to show the rate in the finite time blow-up, the decay condition of initial data on the right is an essential assumption in \cite{MM}.) As the by-product of the method, they were able to obtain estimates on the blow-up rate, which we review below in Section \ref{S:2D-rate}. In \cite{Gong} the authors use the next order approximation of the blow-up profile, see \eqref{E:2d-approx}, it also influences the rate \eqref{E:2d-rate}, which they obtain via the modulation parameters system \eqref{E:2d-mod}. 

As discussed in our paper \cite{KRS1}, investigating the blow-up profiles and rates numerically is a delicate and challenging problem: the numerics can suggest or confirm certain quantities in the analytical models. Thus, in \cite{KRS1} we have confirmed that the profile at the {\it zero order} approximation starts with the ground state \eqref{E:2d-profile} and suggested the initial approximation for the rate, providing numerical examples for the rate $\beta$ (see precise definition in Section \ref{S:2D-rate}) at least as fast as the power $0.5$ (or slightly higher as in the Gaussian example in that paper). The {\it next order} approximation can certainly correct it, as it was shown in \cite{Gong}. 

It is known from, for example, the NLS literature (e.g., \cite{SS1999}), that it is quite challenging to obtain a precise blow-up information numerically in the borderline cases (for example, in the 2D critical (cubic) NLS). 
Determining the blow-up rate numerically requires not only very intricate computational methods and powerful high resolution (to be able to approach sufficiently close to the blow-up time), but also analytical and heuristic steps could be used to attenuate the rate and profile description to further precision.

\subsection{Review of the rate and profile in 2D critical ZK}\label{S:rate-review}
For the two-dimensional cubic ZK equation 
$$
u_t+\partial_{x_1}\big(\Delta_{\mathbb R^2} u+u^3\big)=0,
    \qquad (t,x)\in [0,T)\times \mathbb R^2,
$$
Chen et al in \cite{Gong} described the dynamics of solutions near the
ground-state solitary wave, in particular, in the blow-up regime, they obtained the rate and profile with the {\it next order} approximation; namely, they decompose the solution as
\begin{equation}\label{E:2d-u}
u(t,x)= \frac{1}{\lambda(t)}
\left[Q_{b(t)}(y)+\varepsilon(t,y) \right],
\qquad y=\frac{x-x(t)}{\lambda(t)} ,
\end{equation}
where $Q_b$ is an adjusted profile 
of the form
\begin{equation}\label{E:2d-approx}
Q_b = Q+bP\phi_b,
\end{equation}
with $Q$ being the positive radial solution of the 2D elliptic equation,
\begin{equation}\label{E:2d-profile}
-\Delta Q+Q-Q^3=0 \quad \mbox{in} ~~\mathbb R^2.
\end{equation}
Here, $P$ is the first-order correction (non-localized) profile and $\phi_b$ is a cutoff that localizes the non-decaying part of $P$ (for the details refer to \cite[Sections 1.3, 2]{Gong}). 
The modulation parameters satisfy, at the leading order, the following system
\begin{equation}\label{E:2d-mod}
\frac{\lambda_s}{\lambda}=-b, \qquad b_s+\theta b^2=0, \qquad \frac{ds}{dt}=\frac{1}{\lambda^3},
\end{equation}
where the correction $\theta$ is defined as
\begin{equation}\label{E:2d}    
\theta:= \frac{2\ds \int_{\mathbb R}\frac{|\widehat F(\xi)|^2}{1+|\xi|^2}\,d\xi}{\ds \int_{\mathbb R}|\widehat F(\xi)|^2\,d\xi}, \qquad F(y_2)=\int_{\mathbb R}\Lambda Q(y_1,y_2)\,dy_1,
\end{equation}
recalling that in the 2D cubic ZK equation the scaling generator is $\Lambda Q=Q+y\cdot\nabla Q \equiv Q +r \partial_rQ$.    

Denoting 
\begin{equation}\label{E:theta-comp}
\beta=\frac{1}{3-\theta}, 
\end{equation}
it was shown that in the considered finite-time blow-up regime with negative energy and mass slightly supercritical,  
\begin{equation}\label{E:2d-rate}
\lambda(t)\sim c\,(T-t)^\beta, \qquad \|\nabla u(t)\|_{L^2}
\sim \frac{C}{(T-t)^\beta},
\qquad \beta\in\Big(\frac57,\frac56\Big).
\end{equation}
The authors in \cite{Gong} reported via elementary numerical computations $\theta \approx 1.66032$; then substituting this into \eqref{E:theta-comp} yields $\beta \approx 0.7464$.

In our finite domain computations, we obtain $\theta \approx 1.6608225774$, which then substituted into \eqref{E:theta-comp} gives the (next order) approximation of the blow-up rate $\beta$ as 
$$
\beta \approx 0.747. 
$$

\subsection{Approximation of rate and profile for the 3D critical ZK blow-up}\label{S:3Drate}
Applying the same reasoning as in \cite{Gong} for the 3D critical ZK equation \eqref{ZK}, we use the following near-soliton decomposition 
\begin{equation}\label{E:3DZK-ansatz}
 u(t,x) = \frac{1}{\lambda(t)^{3/2}}
    \left[Q_{b(t)}(y)+\varepsilon(t,y) \right],
    \qquad
    y=\frac{x-x(t)}{\lambda(t)}, 
\end{equation}
where $x,y \in \mathbb R^3$ and $Q_{b}$ is a localized deformation or next order approximation of $Q$ of the form
$$
 Q_{b}=Q + b P \phi_b,
$$
with $Q$ being the positive radial ground state solution of \eqref{E:Q}, i.e., 
$$
-Q +\Delta Q + Q^{7/3}=0 \qquad \text{in }\mathbb R^3.
$$
Analogously, the correction $P$ comes from the linearized
profile equation  
$$
\partial_{y_1}L P=\Lambda Q,
$$ 
where 
$$
L = -\Delta+1-\frac73 Q^{4/3} \quad 
\mbox{and} \quad \Lambda Q=\frac32 Q + y\cdot\nabla Q.
$$
The cutoff function $\phi_b$, similarly 
to the two-dimensional case, localizes the correction $P$ (here, the precise $P$ is not needed in the computation of the blow-up rate).

We can write the corresponding finite-dimensional modulation system in the same form as in \cite{Gong} (which is inspired by the 1D critical gKdV \cite{MMR})
$$
\frac{\lambda_s}{\lambda}=-b,
\qquad b_s+\theta_{3d} b^2=0,
\qquad \frac{ds}{dt}=\frac{1}{\lambda^3}, 
$$
with the constant $\theta_{3d}$ defined as follows: denote the function 
$$
F(y_2,y_3)= \int_{\mathbb R} \Lambda Q (y_1,y_2,y_3)\,dy_1,
$$
then  
\begin{equation}\label{E:theta_3}
\theta_{3d} := \frac{2 \ds \int_{\mathbb R^2} \frac{|\widehat F(\xi_2,\xi_3)|^2}{1+\xi_2^2+\xi_3^2} \, d\xi_2 \, d\xi_3}{\int_{\mathbb R^2} |\widehat F(\xi_2,\xi_3)|^2 \, d\xi_2 \, d\xi_3}.
\end{equation}

To compute this constant in a different way, we take advantage of radiality of the ground state $Q$. We first take the Fourier transform, which
gives
$$
\widehat F(\xi_2,\xi_3)= C \, 
\widehat{\Lambda Q}(0,\xi_2,\xi_3),
$$
with the factor $C=\sqrt{2\pi}$ (not important here, as it cancels in the quotient \eqref{E:theta_3}). Since the ground state $Q$ is
radial, we have 
$$
\Lambda Q =\frac32 Q(r) + r \partial_r Q, \quad r = \sqrt{y_1^2+y_2^2+y_3^2},
$$ 
and thus, we can reduce the computation of $\hat{F}$ via the Fourier projection-slice theorem and the reduction to the radial Fourier transform (e.g., see \cite{GrafakosTeschl2013}) to the one-dimensional radial integral
\begin{align}\label{E:H3}
H_3(k)= & \sqrt{\frac{2}{\pi}} \int_0^\infty \Lambda Q(r) \frac{\sin(kr)}{kr} r^2\,dr,
\qquad k=\sqrt{\xi_2^2+\xi_3^2},\notag \\
= & \frac1{k}  \sqrt{\frac{2}{\pi}} \int_0^\infty \big(\tfrac32 Q(r) +r Q^\prime(r) \big) \sin(kr)\, r \,dr, 
\end{align}
since 
$$ \widehat{F}(\xi_2,\xi_3) = \sqrt{2\pi} H_3(k) \qquad \qquad $$
and 
$$\int_{ {\mathbb S}^2} e^{-i kr \cos \theta} d\omega_{\mathbb S^2} = 
\int_0^{2\pi} \int_0^{\pi} e^{-i kr \cos \theta} \sin \theta \, d\phi d\theta 
= 4\pi \frac{\sin (kr)}{kr}.
$$
Substituting this into \eqref{E:theta_3} and using that $H_3$ is even, we obtain 
\begin{equation}\label{E:theta_3-easy}
\theta_{3d} = \frac{2 \ds \int_0^\infty \frac{|H_3(k)|^2}{1+k^2}\,k\,dk}{\ds \int_0^\infty |H_3(k)|^2\,k\,dk}.     
\end{equation}
The blow-up rate in this case is then given by
$$
\lambda(t)\sim  (T-t)^{\beta}, \qquad \mbox{where} \quad \beta = \frac1{3-\theta_{3d}},
$$
and 
$$
\|\nabla u(t)\|_{L^2(\mathbb R^3)} \sim \frac1{(T-t)^{\beta}},
$$
with $\theta_{3d}$ computed from \eqref{E:theta_3-easy}. 

Using the fact that our ansatz \eqref{E:3DZK-ansatz} has scaling factor $\lambda(t)^{-3/2}$, we obtain the corresponding $L^\infty$-rate 
\begin{equation}\label{E:3D-rate}
\|u(t)\|_{L^\infty(\mathbb R^3)} \sim \frac1{\lambda(t)^{3/2}} \sim \frac1{(T-t)^{\frac32\beta} }.
\end{equation}

In our computations, we get $\theta \approx 1.4481476279$, and thus, 
$$
\beta_{3d} = 0.644, 
$$
(with $L^\infty$-norm rate $\frac32 \beta \approx 0.967$).

\subsection{Radial approach to compute 2D rate correction $\theta$}

We can also employ radiality in 2D critical setting to estimate the 
rate more precisely. 
(Here, $Q$ is the ground state solution to \eqref{E:2d-profile}.)

Recalling that 
$$
F(y_2) = \int_{\mathbb R} \Lambda Q (y_1, y_2) \, dy_1, \quad \mbox{where} \quad \Lambda Q = Q +y \cdot \nabla Q, 
$$
we take the Fourier transform, which gives
$$
\widehat F(\xi_2)= C \, \widehat{\Lambda Q}(0,\xi_2),
$$
with the factor $C=\sqrt{2\pi}$ (which is again not important here, as it cancels in the quotient \eqref{E:2d}). Since the 2D ground state $Q$ is radial, we have 
$$
\Lambda Q = Q(r) + r \partial_r Q, \quad r = \sqrt{y_1^2+y_2^2},
$$ 
and thus, we can reduce the computation of $\hat{F}$ via the Fourier projection-slice theorem and the reduction to the radial Fourier transform (e.g., see \cite{GrafakosTeschl2013}) to the one-dimensional radial integral
$$ 
H_2(k)= \int_0^\infty \Lambda Q(r) \, J_0 (kr) r \,dr, \quad k=|\xi_2|,
$$
since 
$$ 
\widehat{F}(\xi_2) = \sqrt{2\pi} H_2(k)
$$
and 
$$
\int_0^{2\pi} e^{-i kr \cos \theta} d \theta = 2\pi J_0 (kr),
$$
where $J_0$ is the Bessel function of the first kind of order 0.  
Substituting this into \eqref{E:2d}, and also using the fact that $H_2$ is even, we obtain 
\begin{equation}\label{E:theta_2-easy}
\theta_{2d}= \frac{2 \ds \int_0^\infty \frac{|H_2(k)|^2}{1+k^2}\,dk}{\ds \int_0^\infty |H_2(k)|^2\,dk}.
\end{equation}
Using this approach, we also obtain a similar value for $\theta = \theta_{2d}$ as written at the end of \S \ref{S:rate-review}.


\subsection{General dimension formal blow-up analysis}

It is useful to show how the ground state $Q_d$ and the correction $P_d$
enter the blow-up dynamics. 
In dimension $d$, with the
$L^2$-critical exponent
$$
p=1+\frac4d,
$$
we denote by $f_d$ the real homogeneous nonlinearity of
degree $p$ used in the corresponding ZK equation.  
For the dimensions that we are specifically interested here,
$$
f_2(u)=u^3,\qquad
f_3(u)=u^{7/3}=|u|^{4/3}u,
\qquad 
f_4(u)=u^2.
$$
Note that if we wrote a uniform choice $f_d(u)=|u|^{4/d}u$, it would give us
$f_4(u)=|u|u$ in four dimensions, which is not the standard quadratic
nonlinearity $u^2$ (however, the approach works for any combination, since the formal calculation below is the same for any of these choices due to the positive ground state).

We make the change of variables
\begin{equation}\label{E:rescaling-v}
u(t,x)=\frac1{\lambda(t)^{d/2}}v(s,y),
\qquad
y=\frac{x-x(t)}{\lambda(t)},
\qquad
\frac{\dd s}{\dd t}=\frac1{\lambda(t)^3}.
\end{equation}
Substituting into the $d$-dimensional critical ZK equation gives
\begin{equation}\label{E:rescaled-r}
v_s-\frac{\lambda_s}{\lambda}\Lambda_dv
-\frac{x_s}{\lambda}\cdot\nabla v
+\partial_{y_1}\bigl(\Delta v+{f_d(v)}\bigr)=0,
\qquad
\Lambda_dv=\frac d2v+y\cdot\nabla v.
\end{equation}
Choosing in \eqref{E:rescaled-r}, at leading order,
$\frac{(x_1)_s}{\lambda}=1$, and $\frac{(x_j)_s}{\lambda}=0$ for $j\geq 2$, we obtain
\begin{equation}\label{E:rescaled}
v_s-\frac{\lambda_s}{\lambda}\Lambda_dv
+\partial_{y_1}\bigl(\Delta v + {f_d(v)}-v\bigr)=0
\end{equation}
up to lower-order modulation terms. A stationary positive profile therefore
solves
$$
-Q_d+\Delta Q_d+Q_d^{1+4/d}=0,
$$
which gives rise to the ground state as the zeroth-order
approximation of the blow-up profile.
\smallskip

We next define the parameter
$$
b=-\frac{\lambda_s}{\lambda}.
$$
Substituting $Q_d+bP_d$ into the stationary part of
\eqref{E:rescaled}, we see that the terms of order $b$ cancel provided
\begin{equation}\label{E:P-equation}
\partial_{y_1}L_dP_d=\Lambda_dQ_d,
\quad \mbox{where} \quad
L_d=-\Delta+1-\left(1+\frac4d\right)Q_d^{4/d}.
\end{equation}
Therefore, the localized approximation has the form
\begin{equation}\label{E:Qbd}
Q_{b,d}=Q_d+bP_d\phi_b+\cdots .
\end{equation}
Observe that the profile $P_d$ is not fully localized. Indeed, we set
$$
F_d(z)=\int_{\R}\Lambda_dQ_d(y_1,z)\,\dd y_1,
\qquad
z=(y_2,\ldots,y_d),
$$
and assume that the limits as $y_1 \to \pm \infty$ exist, i.e., 
$$
\lim_{y_1\to\pm\infty}P_d(y_1,z):=P_{d,\pm}(z).
$$
Then integration of \eqref{E:P-equation} with respect to $y_1$
gives the relation
$$
(-\Delta_z+1)\bigl(P_{d,+}-P_{d,-}\bigr)=F_d.
$$
We choose the normalization
$$
P_{d,+}(z)=0.
$$
We also take the cutoff $\phi_b$ in \eqref{E:Qbd} to be equal to one in the soliton
core and equal to zero sufficiently far into the non-decaying left tail.
This makes $Q_{b,d}$ localized, while the resulting cutoff terms are supported away from the
blow-up core.

To analyze the blow-up rate, we consider the transverse tail variable $z$ as above and take a function $h_d$ satisfying
\begin{equation}\label{E:hd}
(-\Delta_z+1)h_d=\Delta_zF_d,
\quad \Rightarrow \quad
\widehat h_d(\xi)=-\frac{|\xi|^2}{1+|\xi|^2}\widehat F_d(\xi),
\end{equation}
where $\xi$ is a $(d-1)$-dimensional variable, e.g., can be written as $(\xi_2, ..., \xi_d)$.
With our normalization for $P_d$, we have $P_d(-\infty,z)=-(F_d+h_d)(z)$. Recalling that $L_d \Lambda_d Q_d = -2Q_d$, we obtain 
\begin{equation}\label{E:scalar-identities}
 (P_d,Q_d)=\frac14\int_{\R^{d-1}}|F_d(z)|^2 \, dz.
\end{equation} 
Denoting by $\Psi_b$ the error in the approximate profile, and ignoring the cutoff for the leading order calculation, we have
\begin{equation}\label{E:Psi-b}
\Psi_b = -b \Lambda_d Q_{b,d} + \partial_{y_1} \big( -\Delta Q_{b,d} + Q_{b,d} - f_d(Q_{b,d}) \big).
\end{equation}
Note that \eqref{E:P-equation} cancels terms of order $b$, and hence, $\Psi_b = O(b^2)$ (in the soliton core). Using the leading order representation \eqref{E:Psi-b}, and skew-adjointness of the scaling generator $\Lambda_d$, \eqref{E:P-equation}, and limits of $P_d$, we obtain 
\begin{equation}\label{E:Psi-inner}
(\Psi_b,Q_d)=-\frac{b^2}{2}\int_{\R^{d-1}}\frac{|\widehat F_d(\xi)|^2}{1+|\xi|^2}\,\dd\xi + o(b^2).
\end{equation}
(We mention, that some care would have to be taken for $d>4$, since $1+4/d <2$ in that case; no such concern occurs for $d=2,3,4$.)

Thus, the profile equation contains the balance $b_s P_d - \Psi_b +$ lower-order terms $=0$.
Pairing with $Q_d$ gives
\begin{equation}\label{E:P-Psi}
b_s(P_d,Q_d)-(\Psi_b,Q_d) \approx  0
\end{equation}
(up to higher order than $b^2$ errors), and hence, again up to higher order terms, we have
\begin{equation}\label{E:bd-law}
b_s+\theta_db^2 = 0.
\end{equation}

We concisely explain how the law \eqref{E:bd-law} gives the rate. Integrating \eqref{E:bd-law} together with $\lambda_s/\lambda = -b$ gives, at the leading order,
$$
b(s)\sim \frac1{\theta_d\, s}
\qquad \mbox{and} \qquad
\lambda(s) \sim s^{-1/\theta_d}\quad \mbox{as} \quad s \to \infty.
$$
Since $\dd t = \lambda^3 \dd s$, we get $T-t \sim s^{(\theta_d-3)/\theta_d}$, which is finite precisely when $\theta_d<3$, and  therefore,
\begin{equation}\label{E:beta-d}
\lambda(t) \sim (T-t)^{\beta_d}, \qquad \beta_d = \frac1{3-\theta_d},
\qquad \|\nabla u(t)\|_{L^2} \sim (T-t)^{-\beta_d}.
\end{equation}
Note that $0<\theta_d<2$ for any $d\geq 2$, which can be seen directly from the definition \eqref{E:theta-d} below, since $0< (1+|\xi|^2)^{-1} < 1$, and thus, we obtain that $\beta_d \in \big(\frac13,1\big)$. We remark that in 1D critical gKdV, the transverse space is absent, so the same reduction would give $\theta_1=2, \beta_1 = 1$.


\subsection{General dimension corrections $\theta_d$ and the 4D example}
\begin{figure}[!htb]
\includegraphics[width=.8\textwidth,height=.47\textwidth]{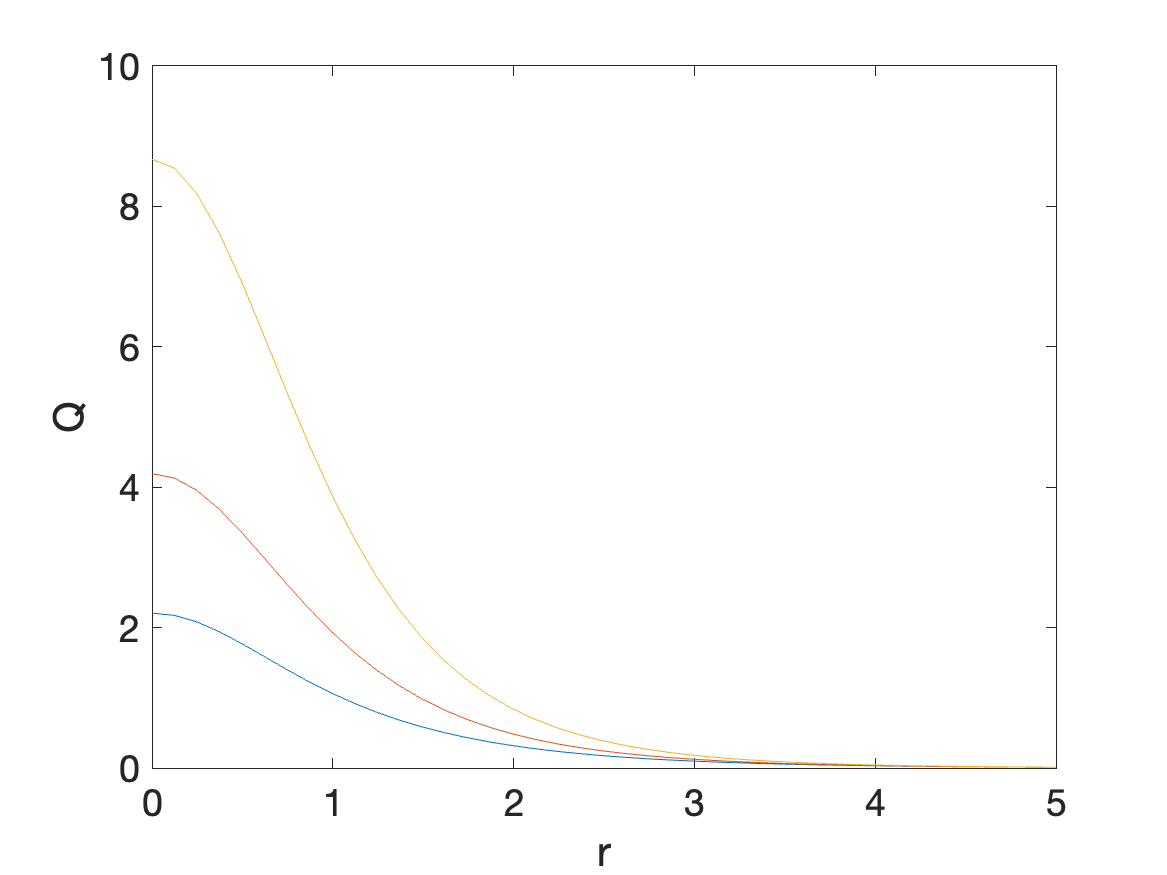} 
\caption{Ground states for the $L^{2}$-critical ZK equation for 
dimensions 4, 3, 2 (from top to bottom).}
\label{NLSsold}
\end{figure}

Using \eqref{E:P-Psi} and \eqref{E:bd-law}, we obtain 
a general formula in dimension $d$ for the correction of the blow-up rate
\begin{equation}\label{E:theta-d}
\theta_{d}= \frac{2 \ds \int_{\mathbb R^{d-1}}
\frac{|\widehat F_d(\xi_2, ..., \xi_d)|^2}{1+\xi_2^2+...+\xi_{d}^2}
        \,d\xi_2\,... d\xi_d}{\int_{\mathbb R^{d-1}}
        |\widehat F_d(\xi_2, ..., \xi_d)|^2 \,d\xi_2\, ... d\xi_d}, 
\end{equation}
where 
$$
F_d (y_2,...,y_d) = \int_{\mathbb R} \Lambda Q (y_1, y_2, ..., y_d) \, dy_1,\qquad 
\Lambda Q = \frac{d}2 Q + y \cdot \nabla Q,
$$
with $Q$ (dropping the subscript $d$) being a ground state solution of the corresponding $d$-dimensional (critical) ground state equation, 
$$
-Q + \Delta_{\mathbb R^d} Q + Q^{\frac4{d} +1} = 0.
$$
One can use the radial property of $Q$ as before to compute $\theta_d$ from \eqref{E:theta-d}.
In particular, in 4D critical (quadratic) ZK, the rate correction $\theta_{4d}$ can be computed as follows (here, $k = \sqrt{\xi_2^2+\xi_3^2+\xi_4^2}$)
\begin{equation}\label{E:theta-4d}
\theta_{4d} = \frac{2 \ds \int_0^\infty \frac{|H_4(k)|^2}{1+k^2}\,k^2\,dk}{\ds \int_0^\infty |H_4(k)|^2\,k^2\,dk},
\end{equation}
where we have used a general formula for radial Fourier transform \cite{GrafakosTeschl2013} (with the Bessel function $J_{\frac{d}2-1}$ of the first kind) 
$$
\hat{f}(k) = k^{1-\frac{d}2} \int_0^\infty f(r) J_{\frac{d}2-1}(kr) r^{\frac{d}2} \, dr,
$$
and thus, in 4D we have
$$
H_4(k) = \frac1{k} \int_0^\infty \Lambda Q (r) \, J_1(kr) \, r^2 \, dr,
$$
where in this case 
$$
\Lambda Q (r) = 2 Q (r) + r \partial_r Q,
$$ 
with $Q$ being the ground state (smooth radial positive decaying 
towards infinity) solution of 
$$
-Q+\Delta_{\mathbb R^4} Q + Q^2 \equiv Q_{rr}+\frac3{r}Q_r -Q +Q^2=0.
$$ 
We show the ground states for various dimensions in the $L^{2}$ critical case in Figure ~\ref{NLSsold}.
\smallskip

For $d=4$, using either \eqref{E:theta-4d} or \eqref{E:theta-d}, we obtain $\theta_{4d} \approx 1.266$, and thus, $\beta_{4d} = \frac1{3-\theta_{4d}} \approx 0.577$.

\section{Properties of ground states in 3D ZK}\label{S:ground}
Recalling the ground state equation for a general power $p$ in 3D, 
\begin{equation}\label{E:Qgen}
\Delta_{\mathbb R^3} Q - Q  + Q^{p} = 0,
\end{equation}
and using the fact that $Q$ is radial, we write $Q(x) = Q(r)$, $r=|X|$, $X \in \mathbb R^3$ and rewrite the above equation as 
$$
Q^{\prime\prime}(r) + \frac2{r} Q^\prime(r) - Q(r) +Q(r)^p = 0,
$$
with $Q'(0) = 0, Q(r) \to 0$ as $r \to \infty$. We also have $Q(r)>0$ (positivity) and monotone decay $Q'(r)<0$ for $r>0$. 

At infinity, due to the decay, the $Q^p$ becomes negligible, so a valid approximation would be $\Delta Q - Q \approx 0$, which gives the decay 
$Q(r) = c \, \frac{e^{-r}}{r} (1+O(\frac1{r})) \quad \mbox{and} \quad Q'(r) = -Q(r) (1+\frac1{r}+O(\frac1{r^2}))$. 

Using this asymptotics, for numerical purposes on a large truncated interval $[0,R]$ a possible boundary condition is 
$$
Q'(R) +(1+\frac1{R})Q(R) \approx 0.
$$

We next record the usual Pokhozhaev identities (either multiplying \eqref{E:Qgen} by $Q$ or by $X \cdot \nabla Q$):
$$
\|Q\|^2_{L^2} + \|\nabla Q\|^2_{L^2} = \|Q\|^{p+1}_{L^{p+1}}
$$
and
$$
\frac12 \|Q\|^2_{L^2} + \frac16 \|\nabla Q\|^2_{L^2} = \frac1{p+1} \|Q\|^{p+1}_{L^{p+1}}.
$$
Solving for $\|\nabla Q\|^2_{L^2}$ and $\|Q\|^{p+1}_{L^{p+1}}$, we obtain
\begin{align}\label{E:P1}
\|\nabla Q\|^2_{L^2} &= \frac{3(p-1)}{5-p}\|Q\|^{2}_{L^{2}},\\
\|Q\|^{p+1}_{L^{p+1}} &= \frac{2(p+1)}{5-p} \|Q\|^2_{L^2} \label{E:P2}
\end{align}
Thus, with the energy convention \eqref{EC}, we obtain 
$$
E[Q] = \frac{3p-7}{2(5-p)}M[Q].
$$
In the $L^2$-critical case $p=\frac73$, we have $E[Q]=0$, which we confirm  numerically as well. We also deduce 
\begin{equation}\label{E:Grad-Pot}
\|\nabla Q\|_{L^2}^2 = \frac32 M[Q] \quad \mbox{and} \quad 
\|Q\|^{10/3}_{L^{10/3}} = \frac52 M[Q].
\end{equation}
In this critical case, we numerically obtain 
\begin{equation}\label{E:L2}
\|Q\|_{L^2} \approx 7.9864, \qquad M[Q] \approx 63.78, \qquad Q(0) \approx 4.1917. 
\end{equation}

We also mention that the ground state $Q$ in this critical case provides a threshold condition for the global existence, which comes from $Q$ being the optimizer for the 3D sharp Gagliardo--Nirenberg inequality ~\cite{W83}
\begin{equation}\label{E:GN}
\|u\|^{10/3}_{L^{10/3}(\mathbb R^3)}
\leq C_{\mathrm{GN}} \|u \|^{4/3}_{L^2(\mathbb R^3)} \|\nabla u\|^{2}_{L^2(\mathbb R^3)}.
\end{equation}
Using \eqref{E:P1}-\eqref{E:P2} from Pokhozhaev identities, the sharp constant can be written as
\begin{equation}\label{E:GN-constant}
 C_{\mathrm{GN}} =\frac53\, M[Q]^{-2/3} \equiv \frac53\,\|Q\|^{-4/3}_{L^2}.
\end{equation}
Bounding the potential term in the energy from \eqref{E:GN} and using the sharp constant \eqref{E:GN-constant}, we get 
\begin{equation}\label{E:coercive-lower}
 E[u]\ge\frac12\left[1-\left(\frac{M[u]}{M[Q]}\right)^{2/3}\right]
 \|\nabla u\|_{L^2}^2,
\end{equation}
and hence, if $M[u]<M[Q]$, we obtain the bound on the gradient
\begin{equation}\label{E:energy-coercivity}
 \|\nabla u\|^{2}_{L^2(\mathbb R^3)} \leq 2 E[u]
\left[1-\bigg(\frac{M[u]}{M[Q]}\bigg)^{2/3}\right]^{-1}.
\end{equation}
Combined with the appropriate local theory, one would get that  
the solutions are global in $H^1(\mathbb R^3)$, under this mass threshold.

With that we describe our numerical approach and then show the results of our simulations.

\section{Numerical approach}\label{S:Num}

Scientific computing in three dimensions poses substantial memory and 
significant computational  challenges. We recall that the solitary 
wave profiles in the ZK equation (in particular, in 3D) are 
exponentially decreasing, and hence, we can treat them using Fourier 
methods. Indeed, on a sufficiently large domain functions that 
decrease fast towards infinity can be continued periodically as 
numerically smooth (that is, to a finite precision). We thus compute on a large $3$-torus. The critical phenomena require high spatial resolutions, and the fact that we are working in three dimensions makes the total number of points grow fast with even a small increase in resolution.  We have found that larger resolution is necessary in the principal direction $x$, so we work with an $xyz$ grid of 
$$
N_x \times N_y \times N_z = 2^{12}\times2^9\times2^9,
$$ 
which is slightly more than one billion grid points. The size of the domain is $8\pi \times 8\pi \times 8\pi$. 
Spectral methods, including Fast Fourier transforms (FFT), are global, that is, they use 
all-to-all mapping for each transformation, and the principal 
bottleneck is in passing information. Using standard computational 
systems or clusters is far from optimal, as the main time cost is in 
communication between nodes. To overcome this, we use a \emph{dense} 
computer system, consisting of several graphical processing units 
(GPUs) on the same node, connected to each other through a so-called 
bridge, which allows for passing information directly between the 
GPUs at high speed. 

Let
$$
\widehat f(k)=(2\pi)^{-3/2}\int_{\R^3}e^{-iX\cdot k}f(X)\,\dd X,
 \qquad k=(k_x,k_y,k_z).
$$
With this convention, we treat the 3D ZK equation \eqref{ZK} in the 
Fourier domain as follows
\begin{equation}\label{ZK-fourier}
 \widehat u_t  =ik_x[(k_x^2+k_y^2+k_z^2)\widehat u  -\widehat{|u|^{4/3}u}].
\end{equation}
The nonlinearity is computed as the FFT of $|u|^{7/3}\mbox{sign}(u)$ 
to avoid problems with the branching of the root due to rounding 
errors. 

Due to the linear term, this is a \emph{stiff} system which loosely 
speaking means that explicit methods are inefficient because of 
stability conditions. Since the stiffness is in the linear part, 
there are many efficient integrators, see for instance  \cite{etna} 
and references therein. For KdV-type equations, exponential time 
differencing proved to be the most efficient, but only if so-called 
$\phi$-functions are precomputed. In the present case this is not 
recommended, since this would require the storage of several 3D 
functions being memory consuming. Therefore, we 
use here an integrating factor method: this means we write \eqref{ZK-fourier} in the form
\begin{equation}	\label{ZKIF}
\frac{d}{dt}	( e^{-ik_{x}(k_{x}^{2}+k_{y}^{2}+k_{z}^{2})t} \, \hat{u}(k,t)) = -ik_{x} \,e^{-ik_{x}(k_{x}^{2}+k_{y}^{2}+k_{z}^{2})t}\, \widehat{u^{7/3}}(k,t),
\end{equation}
to which we apply a standard Runge-Kutta fourth order method (RK4). 
This choice is made in order to make most use of the available GPUs 
memory and maximize the resolution. Compared with exponential time-differencing formulas that store several three-dimensional $\phi$-functions, the integrating-factor formulation reduces GPU memory use while treating the stiff linear part exactly.
Furthermore, this is an explicit method and we do not need to iterate.  The integrating factor lowers 
the stability requirements, but we still need very small time-steps to keep the method stable. 

GPUs mass parallel capacity offers about two orders of magnitude speed-up to standard scripted code (e.g., Matlab or Python), however, direct comparison is not at all trivial, especially, at high resolutions, as the code structure differs significantly. This is only complicated by the actual hardware on which the codes are ran: we use a system with two Intel Xeon E5 processors and two NVIDIA Quadro RTX 8000 GPUs with a NVLink bridge between them. 

We also mention that the code in this paper was tested against the code 
used in \cite{KRS_cyl} to find solutions with radial symmetry in the 
$yz$ plane by the means of Chebyshev collocation. In the relevant 
cylindrically symmetric cases the results agree to $10^{-8}$ and even better if we are away from the blow-up.

\section{Numerical results}\label{S:results}

\subsection{Perturbations of the ground state}
We first consider initial data of type 
$$
u_0(X) = A\, Q(X), \quad A \approx 1, \quad X = (x,y,z).
$$

Taking $A<1$, the simulations show decay of the localized bump 
resolving into radiation, similar to what we showed in the 
cylindrically symmetric computations \cite{KRS_cyl}. Thus, we 
consider $A>1$. For that we take initial data $u_0 = 1.1 Q$. This 
choice is dictated by the observation that the closer the initial 
data are to a soliton, i.e., the factor $A \searrow 1$, the longer it takes for the dynamics to resolve, which can be seen from the scaling property. This makes it increasingly difficult to reach the asymptotic regime. For the case of $A=1.1$, the solution blows up around $t = 4.6$.  

We use the conservation of mass $\Delta_M= \log\big( |M[u_0] - M[u(t)|/M[u_0] \big)$ to control the precision of the temporal evolution scheme. Once the value $\Delta_M$ gets over $-2$, we stop the code. Generally, this quantity overestimates the global error by several orders of magnitude. This way we are able to track this solution until $t=4.5$ (which can affect the approximation of the obtained rates and profiles). 
\begin{figure}[!htb]
\includegraphics[width=0.3\hsize]{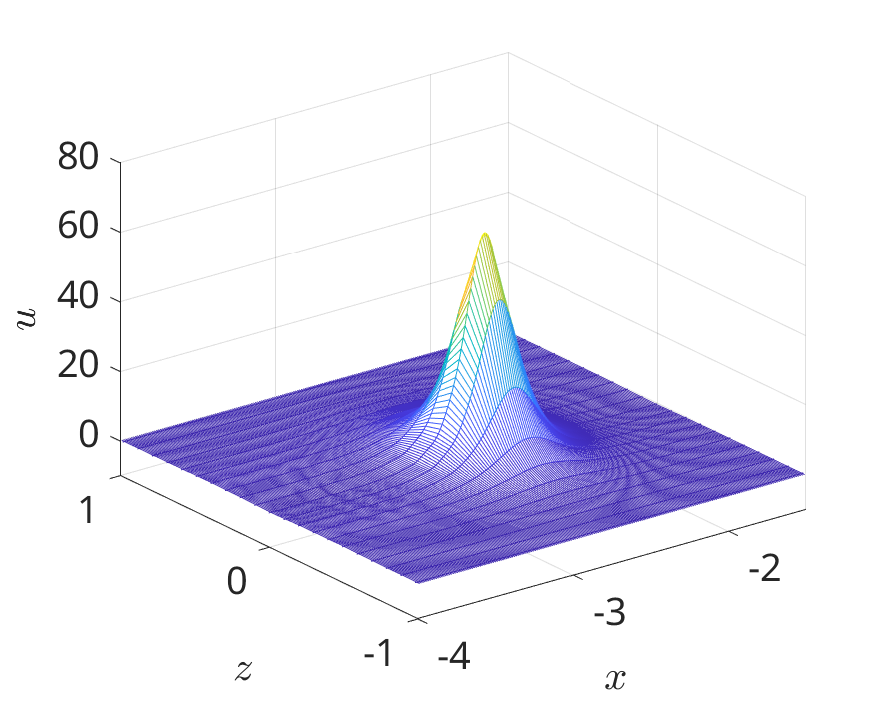}
 \includegraphics[width=0.3\hsize]{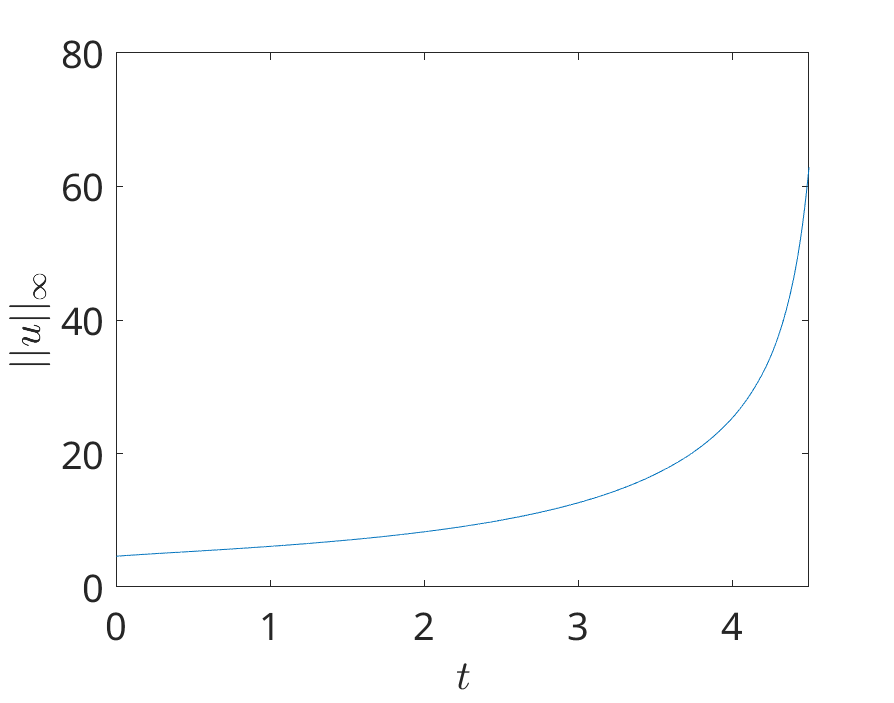}
\includegraphics[width=0.3\hsize]{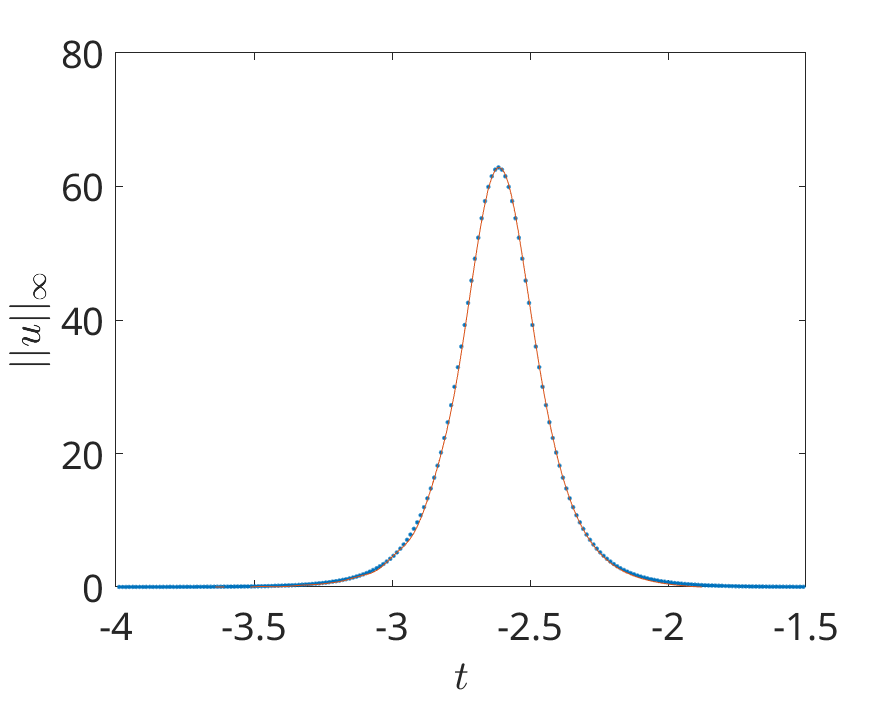}
\caption{Evolution of the perturbed-soliton initial data $u_0 = 1.1Q$. 
On the left a zoom-in of the solution at $t = 4.5$; in the middle the evolution of $\|u(t)\|_{L^\infty}$; on the right a match of the solution at $t = 4.5$ (red line) and a rescaled soliton (blue dots), providing evidence for the asymptotic profile.}
\label{F:Q-1}
\end{figure}

\begin{figure}[!htb]
\includegraphics[width=0.32\hsize]{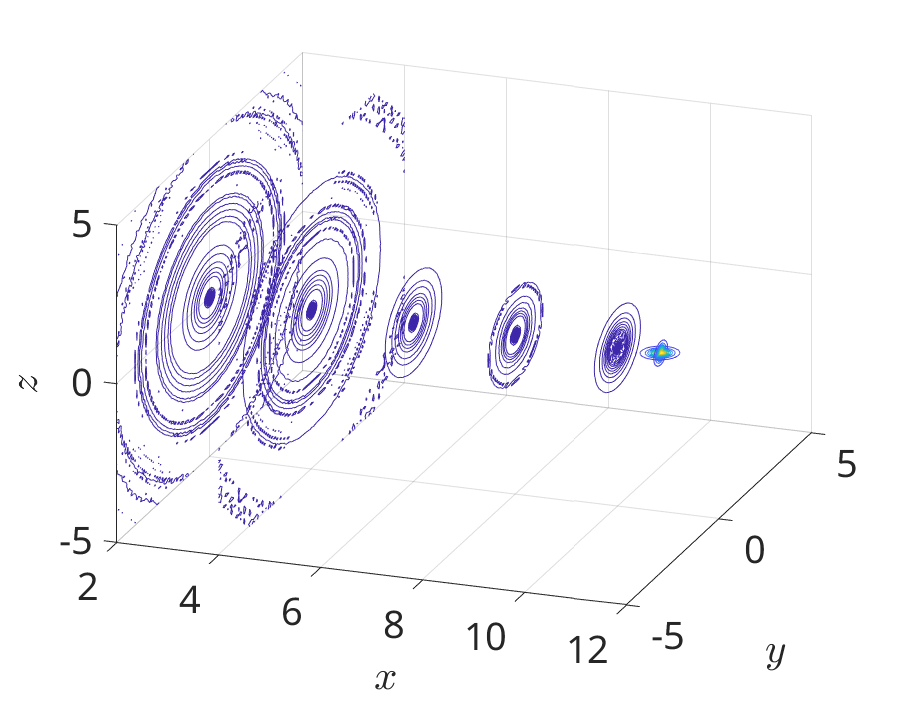}
\includegraphics[width=0.32\hsize]{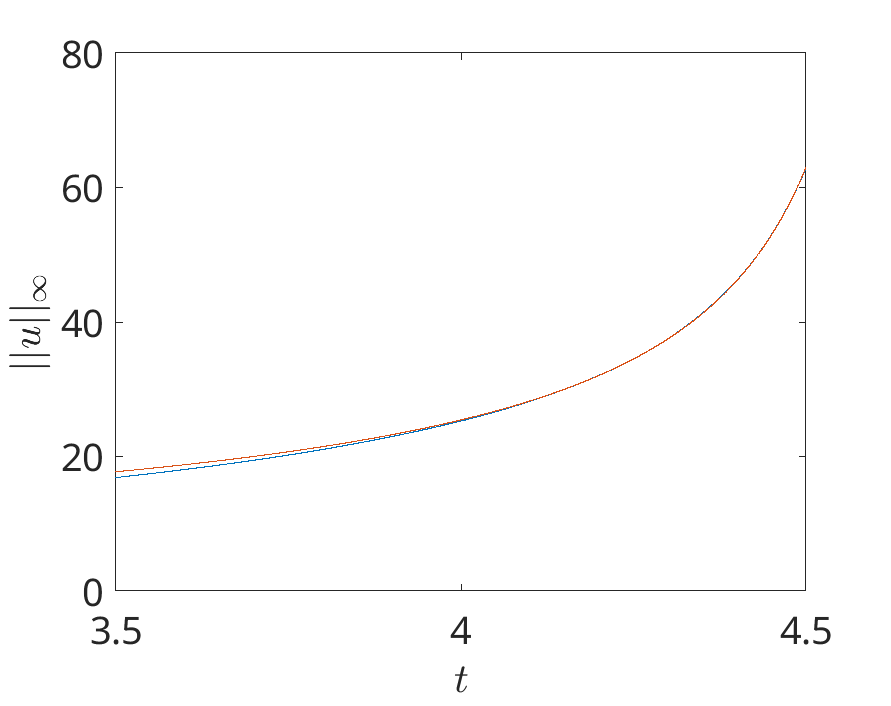}
\includegraphics[width=0.32\hsize]{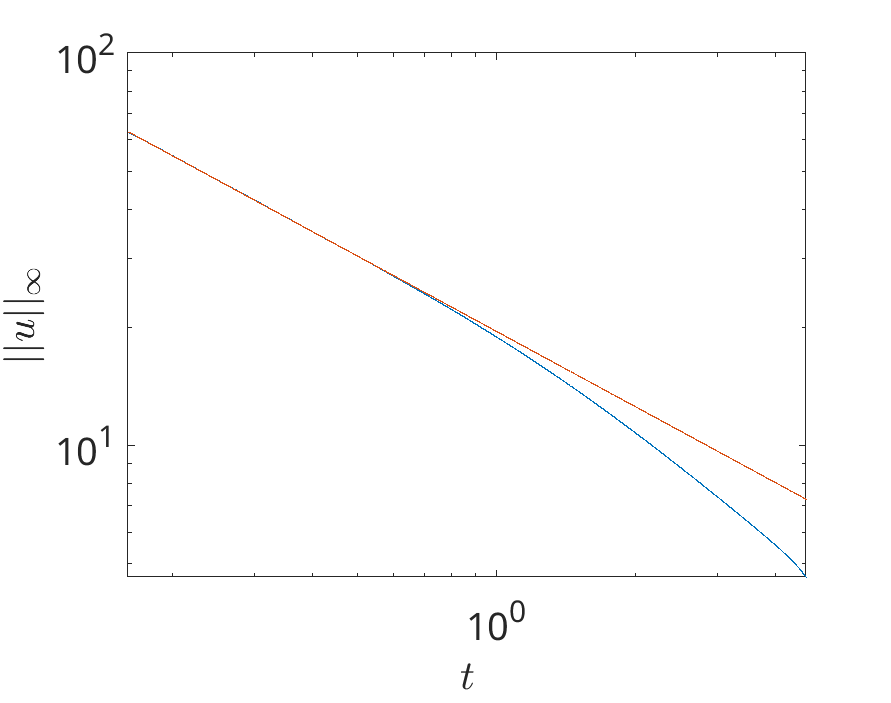}
\caption{On the left a contour plot for the solution with inital data $u = 1.1Q$ at $t = 4.5$. The point of view is rotated to $90^o$ to improve presentation. We have also taken advantage of the periodicity in $x$ to the same end. In the middle the $\|u(t)\|_{L^\infty}$ norm in blue and a blowup rate following $b(t_c-t)^a$ in red. On the right the same two curves in a log-log scale. We found $t_c = 4.66$, $a = -0.64$, $b = 19.5$. }
\label{F:Q-2}
\end{figure}

We show some time evolution and tracking of the $L^\infty$ norm as well as the matching with the rescaled ground state profile in Figure ~ \ref{F:Q-1}. In the next Figure \ref{F:Q-2} we provide a contour plot of the solution at $t=4.5$, noticing the higher concentration (see in blue) of the solution as it moves along the $x$-axis. On the same figure we also show the time dependence of the $L^\infty$ norm, which is connected with the blow-up rate, and a comparison with the power law $b (t_c-t)^a$. Since the final computational time $t=4.5$ is not yet sufficiently close to the fitted blow-up time $t_x \approx 4.66$, the fitted exponent $a$ shows the pre-asymptotic regime.

\subsection{Gaussian initial data}\label{S:Stability}

We next consider 
\begin{equation}\label{ID:Gauss}
u(x,y,z,0) = A \, e^{-\alpha(x^2+y^2+z^2)}, ~~ A \gg 1.
\end{equation}
We note that 
$$
M[u] = A^2 \, \big(\frac{\pi}{2\alpha} \big)^{3/2}.
$$ 
Hence, the mass-threshold amplitude is 
$$
A_M(\alpha) =  \|Q\|_{L^2} \big(\frac{2\alpha}{\pi} \big)^{3/4}.
$$
Then $M[u_0] < M[Q]$ when $A<A_M(\alpha)$. For $\alpha=1$, we get $A_M(1) \approx 5.69$. 
Thus, the dispersion regime for Gaussian initial data with $\alpha=1$ is 
\begin{equation} 
A < A_M(1) \approx 5.69.
\end{equation}

We also compute the energy (to check when it is negative):
\begin{equation}\label{E:Gaussian-energy}
E[u_0]=\frac32\alpha A^2\left(\frac\pi{2\alpha}\right)^{3/2}
 -\frac3{10}A^{10/3}\left(\frac{3\pi}{10\alpha}\right)^{3/2}.
\end{equation}
The energy changes sign at
\begin{equation}\label{E:Gaussian-energy-threshold}
A_E(\alpha)=\frac59\,15^{7/8}\alpha^{3/4}, \qquad A_E(1)\approx 5.94.
\end{equation}
Hence, we expect solutions with Gaussian initial data, $\alpha=1$, and $A> A_E(1) \approx 5.94$ to blow-up in finite time. 

Thus, the choice that we use in our simulations, $A=6.5$, $\alpha=1$, has both supercritical mass and negative energy,
$$
u_0=6.5 e^{-(x^2 + y^2 + z^2)}.
$$ 
This initial condition produces the time evolution that shows the concentration, indicating the blow-up at a finite time, at approximately $t^* =0.85$. We are able to run the simulations until $t=0.8$ with $16 000$ time-steps. 
\begin{figure}[!htb]
\includegraphics[width=0.32\hsize]{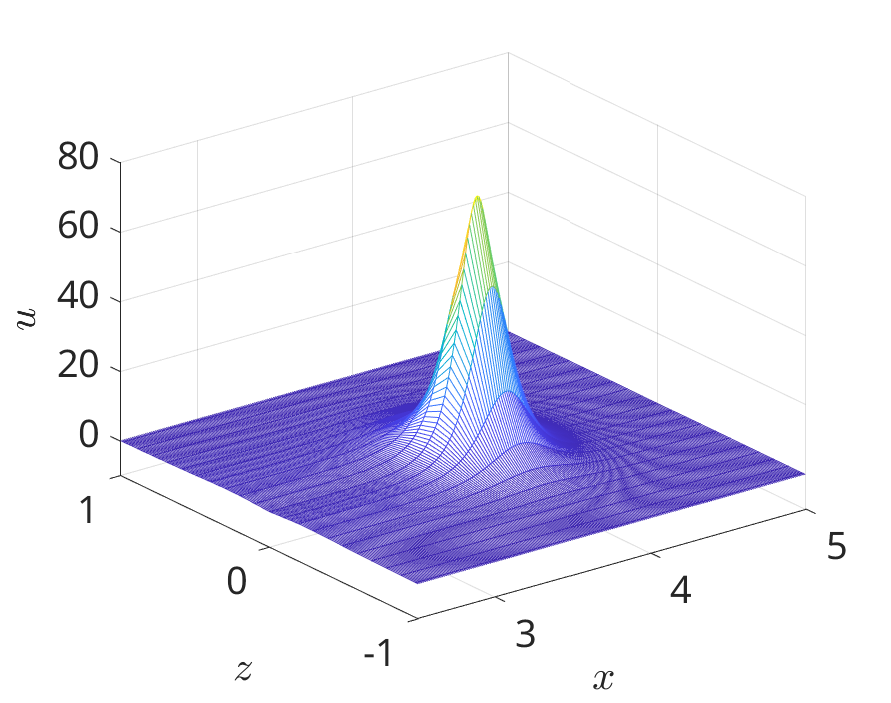}
 \includegraphics[width=0.32\hsize]{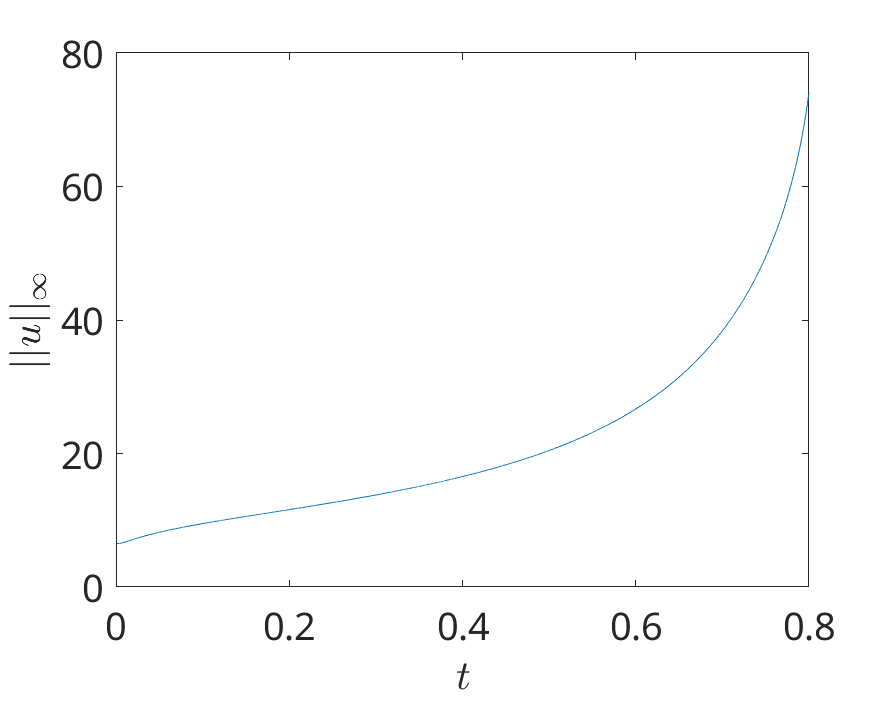}
\includegraphics[width=0.32\hsize]{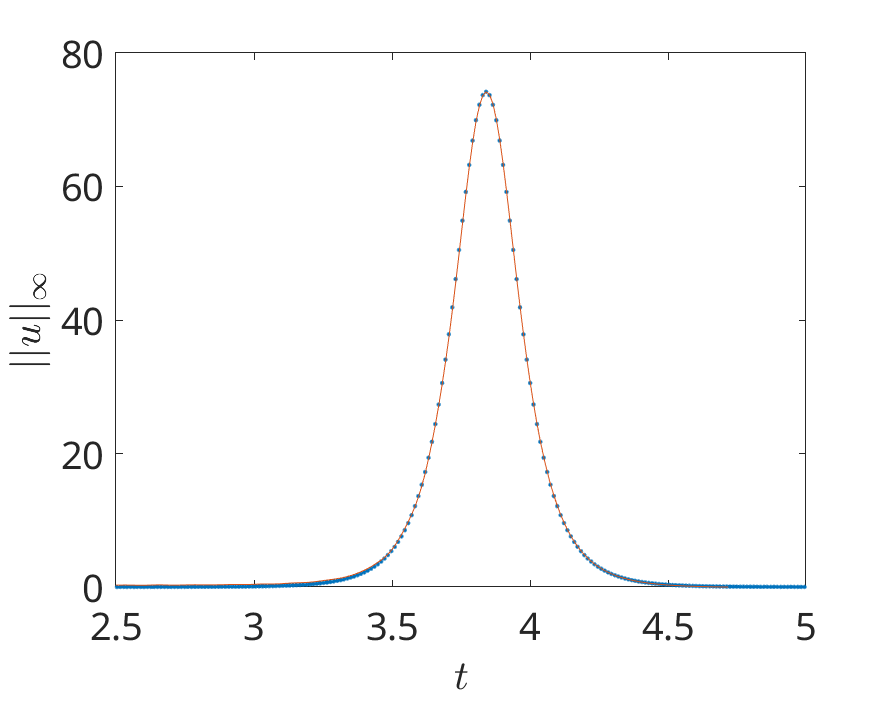}
\caption{Evolution of Gaussian initial data $u_0 = 6.5e^{-(x^2+y^2+z^2)}$. On the left a zoom-in of the solution at $t = 0.8$; in the middle the evolution of $\|u(t)\|_{L^\infty}$; on the right a match of the solution at $t = 0.8$ (red line) and a rescaled soliton (blue dots), providing evidence for the asymptotic profile.}
\label{F:G-1}
\end{figure}
In Figure ~\ref{F:G-1} we show the evolution of the Gaussian initial condition, its $L^\infty$ norm in time and the profile at the final computational time with the matching of the rescaled ground state to show a possible asymptotic profile.

\begin{figure}[!htb]
\includegraphics[width=0.32\hsize]{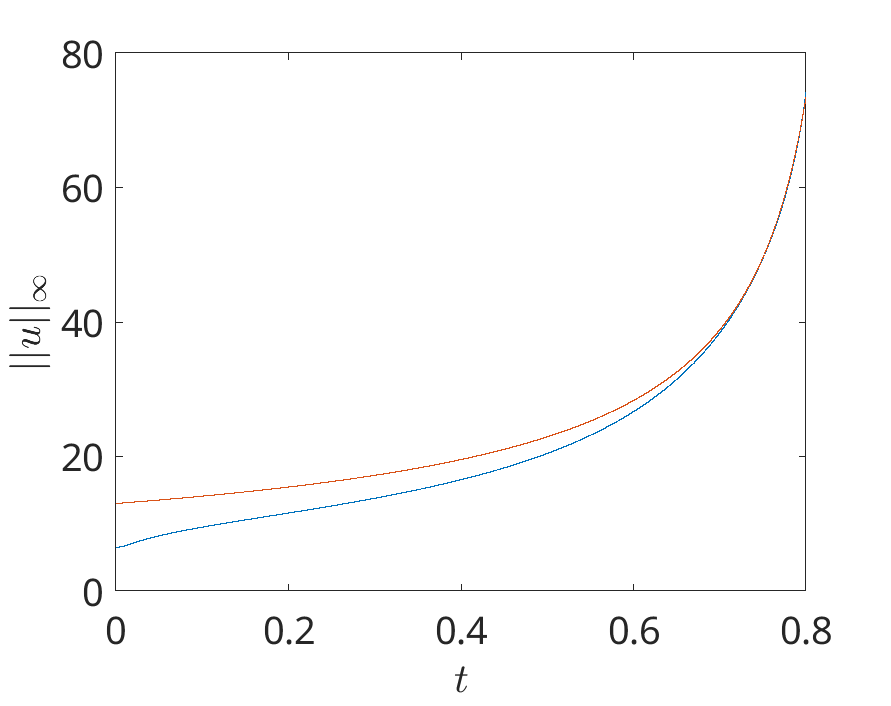}
\includegraphics[width=0.32\hsize]{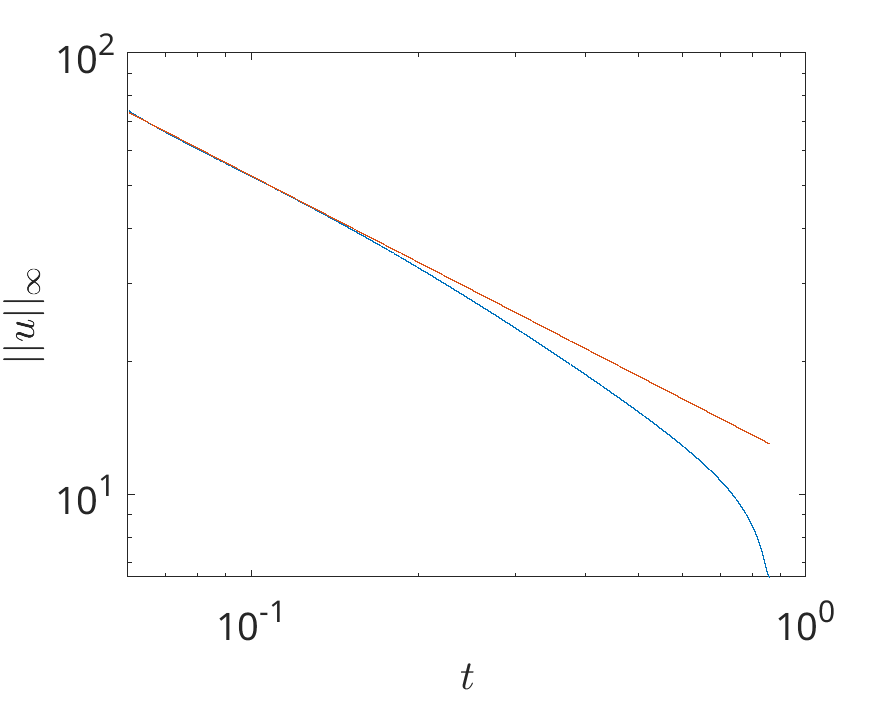}
\caption{On the left $\|u(t)\|_{L^\infty}$ norm in blue and a blowup rate following $b(t_c-t)^a$ in red. On the right the same two curves in a log-log plot, we found that $t_c = .86$, $a = -0.65$, $b = 11.8$.}
\label{F:G-2}
\end{figure}
In Figure ~\ref{F:G-2} we compare the blow-up rate with the power law $b (t_c-t)^a$ and approximate the relevant constants. As it was mentioned, since in the current computations we are reaching only the time $t=0.8$ with the blow-up time around $t_c=0.86$, we are still in a pre-asymptotic regime, which can be seen in a difference when comparing with rate $\beta_{3d}$ (and the corresponding quantity for the $L^\infty$ rate).

\subsection{Double Gaussian initial data}\label{S:doubleG}
Next, we are interested in the setting that breaks the $yz$-symmetry. To this end we take two-bump initial condition, 
\begin{equation}\label{E:doubleG}
u_0 = 5 e^{-(x^2 + y^2 +(z-\pi/2)^2)} + 5 e^{-(x^2 + y^2 +(z+\pi/2)^2)}.
\end{equation}
These are two overlapping Gaussians, with the peaks separated by $\pi$ in the $z$-direction. Each individual Gaussian has a subcritical mass, 
$$
25\left(\frac\pi2\right)^{3/2}\approx 49.22 \approx 0.77\, M[Q].
$$
Including the positive overlap term, the total mass is
\begin{equation}\label{E:doubleG-mass}
M[u_0]=50\left(\frac\pi2\right)^{3/2} \big(1+e^{-\pi^2/2}\big)
 \approx 99.14 \approx1.55 \, M[Q].
\end{equation}
Thus, both bumps are individually subcritical, while their superposition is well above the ground-state mass. We also compute the energy of this double Gaussian,
$$
E[u_0] \approx 27.4685 >0,
$$
noting that if we treated two bumps separately, then each bump would have $E[5e^{-X^2}]\approx 15.15$, thus, the overlap takes only a small amount of energy, since Gaussians decay very fast (compare this with the energy in the next example).

\begin{figure}[!htb]
\includegraphics[width=0.32\hsize]{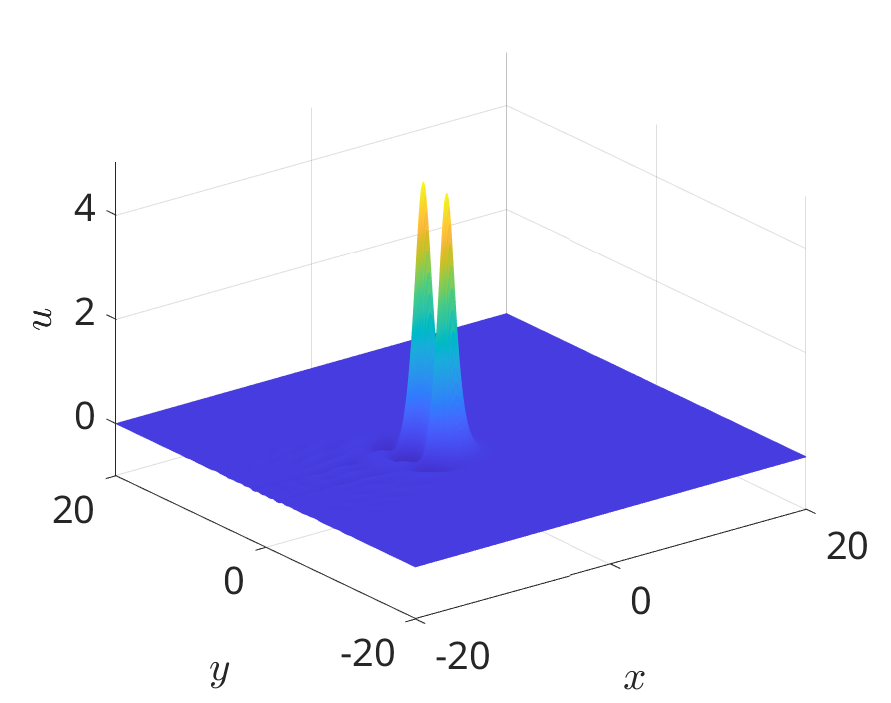}
 \includegraphics[width=0.32\hsize]{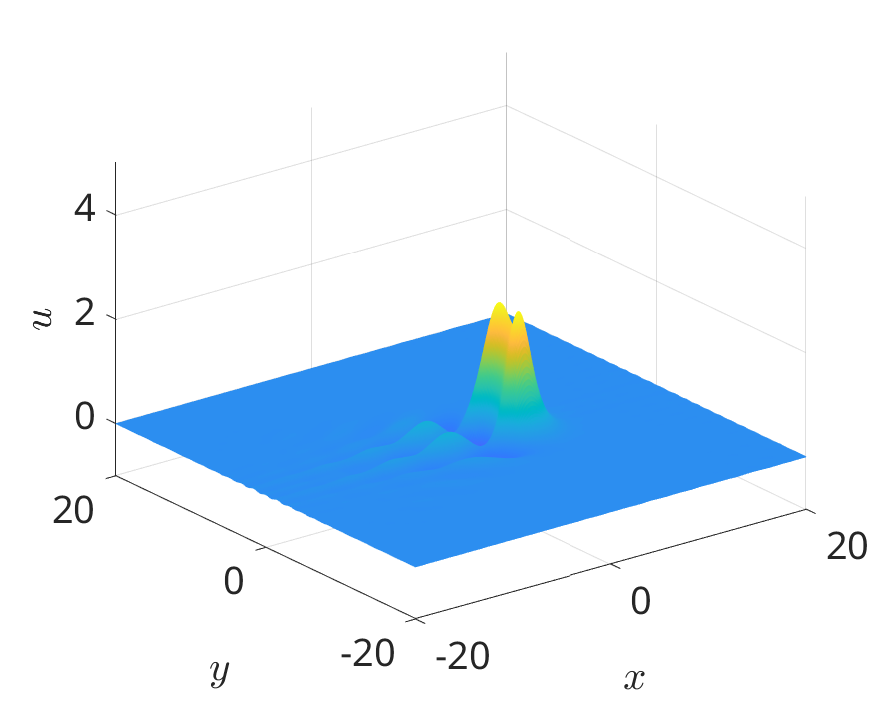}
  \includegraphics[width=0.32\hsize]{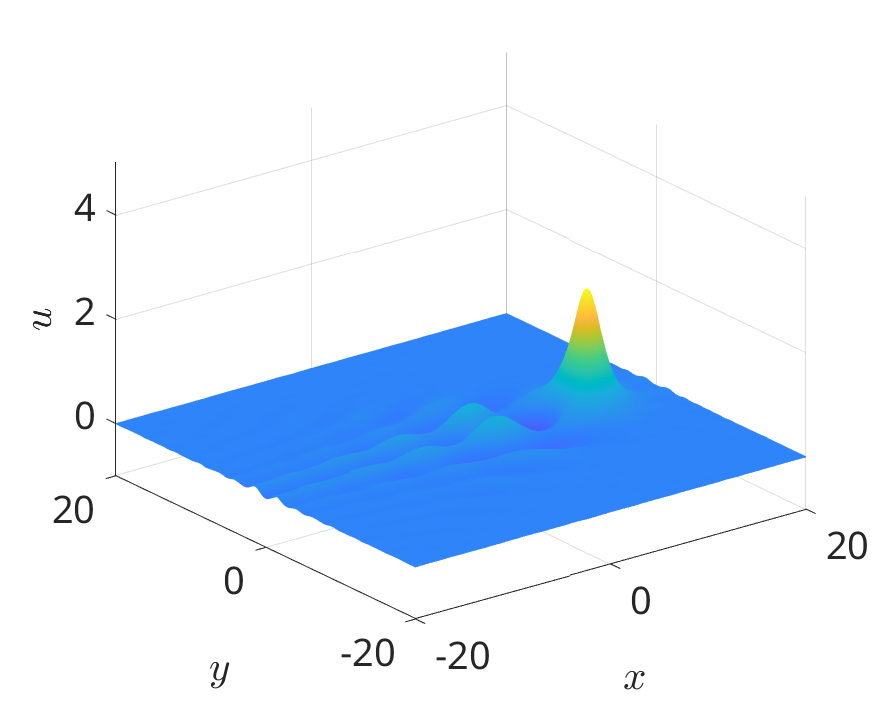}
\caption{Evolution of the double Gaussian initial condition \eqref{E:doubleG} at $t=0.5$ on the left, $t=2$ in the middle, $t= 5$ on the right.}
\label{F:doubleG-1}
\end{figure}

\begin{figure}[!htb]
\includegraphics[width=0.32\hsize]{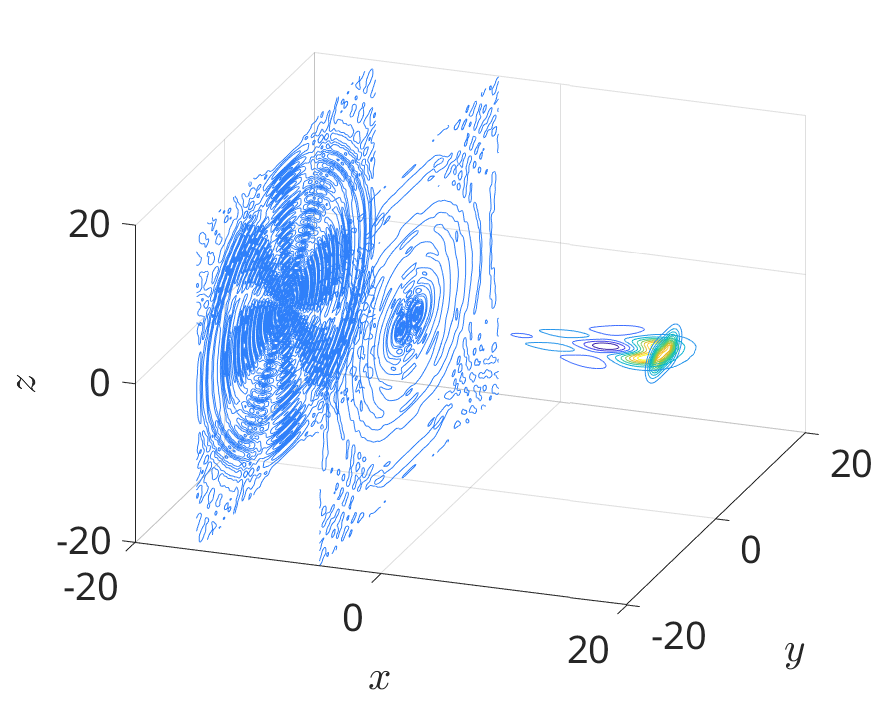}
 \includegraphics[width=0.32\hsize]{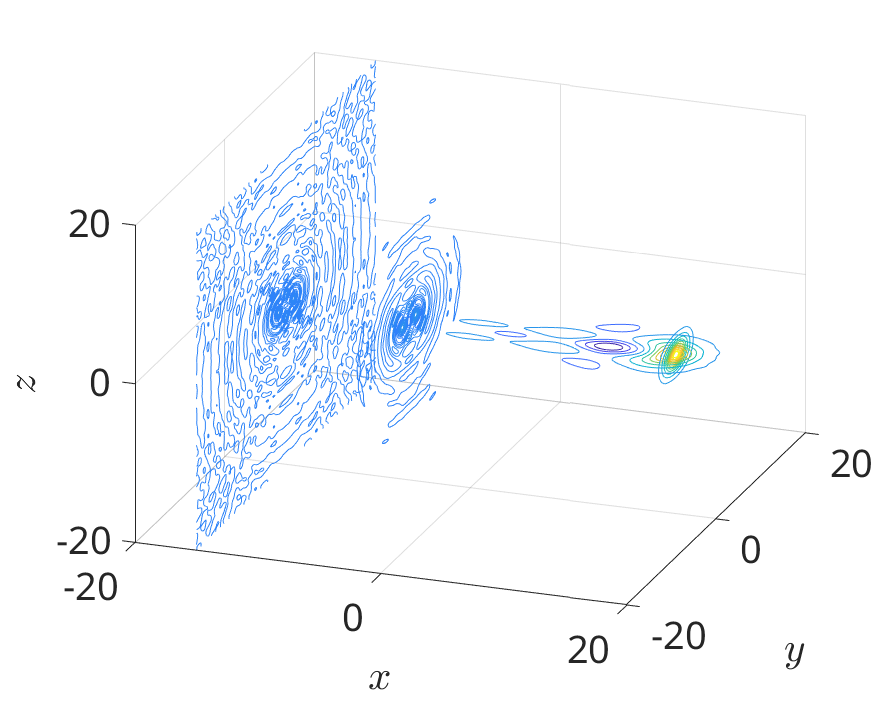}
  \includegraphics[width=0.32\hsize]{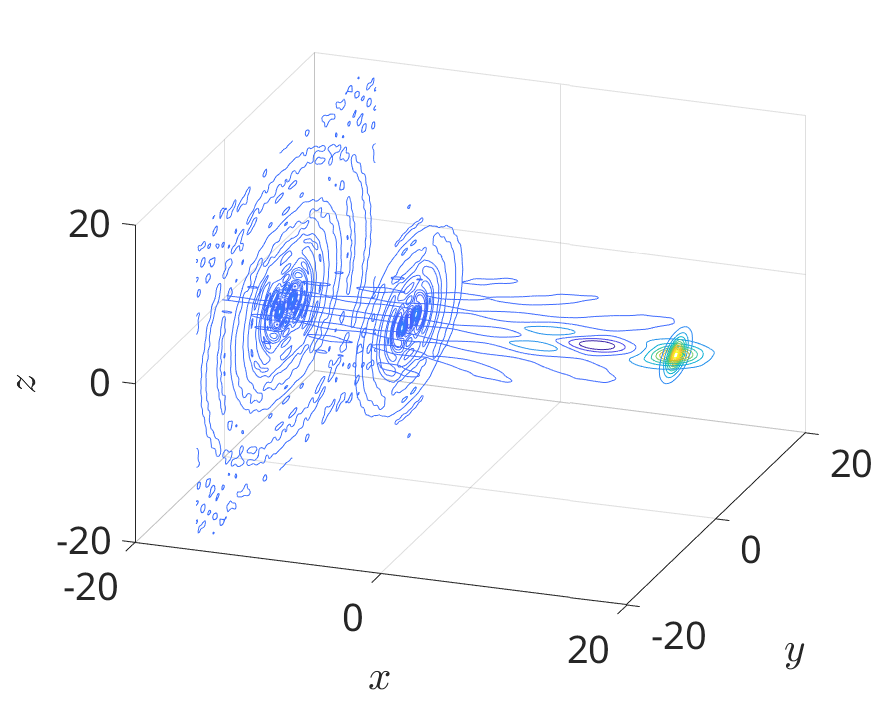}
\caption{Contour plots for the evolution of double Gaussian from \eqref{E:doubleG} at $t=2$ on the left, $t=3.5$ in the middle, $t= 5$ on the right. The drawings are done in the planes $z =0$, $x = -15$, $x = -5$ and the maximum of $u$. There are 10 isolines on each plane; note that the direction of the $x$-axis is reversed compared to the standard plots to improve presentation. We use the periodicity in $x$ as explained in the text.}
\label{F:doubleG-2}
\end{figure}

The time evolution can be seen in Figure ~\ref{F:doubleG-1} and the corresponding contour plots in Figure~ \ref{F:doubleG-2}.


\begin{figure}[!htb]
\includegraphics[width=0.38\hsize,height=0.27\hsize]{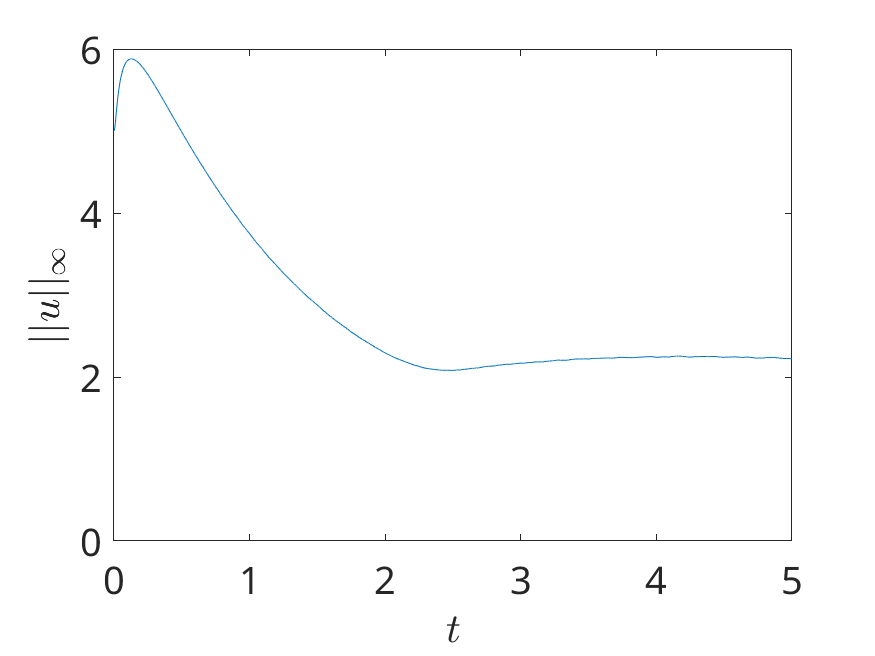}
\caption{Evolution of the $\|u(t)\|_{L^\infty}$ for the double Gaussian initial condition \eqref{E:doubleG}.}
\label{F:doubleG-3}
\end{figure}

Our simulations show that as time progresses, two Gaussians merge into a single localized peak, after which the maximum decreases and the solution thermalizes, approaching a bounded state, see the last plot in Figure~\ref{F:doubleG-1} and also the time dependence of the $L^\infty$ norm in Figure~\ref{F:doubleG-3}.    
This example shows that total mass above $M[Q]$ is not a sufficient numerical criterion for blow-up.

\subsection{Double soliton initial data}\label{S:doubleQ}
\begin{figure}[!htb]
\includegraphics[width=0.32\hsize]{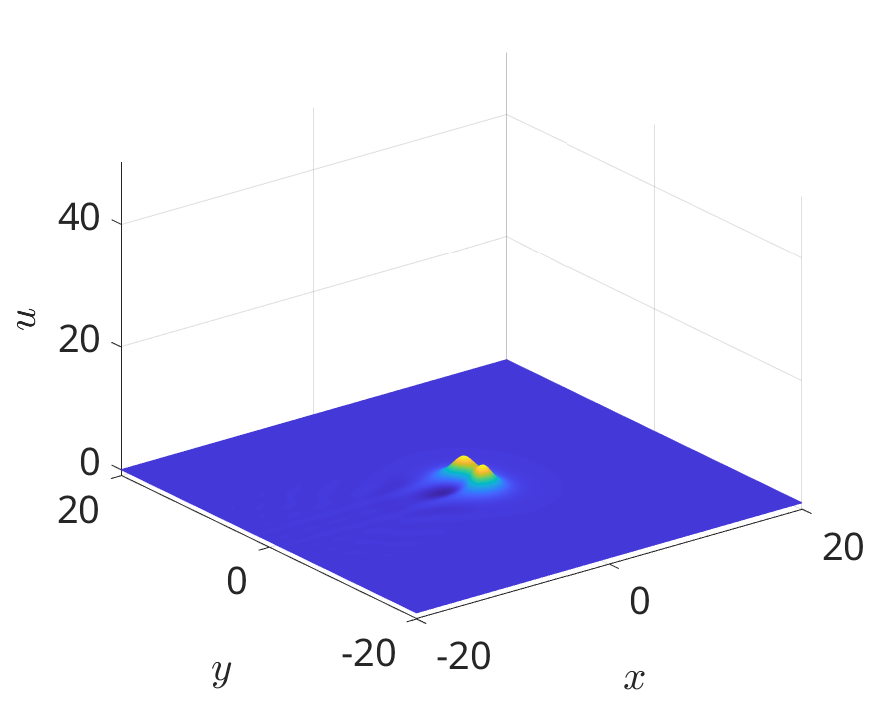}
 \includegraphics[width=0.32\hsize]{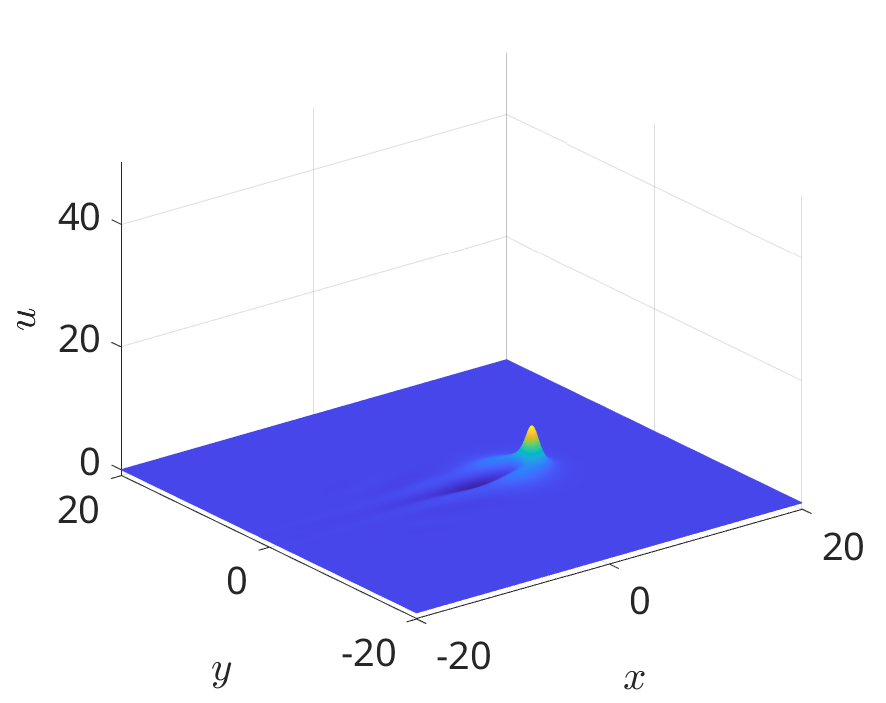}
  \includegraphics[width=0.32\hsize]{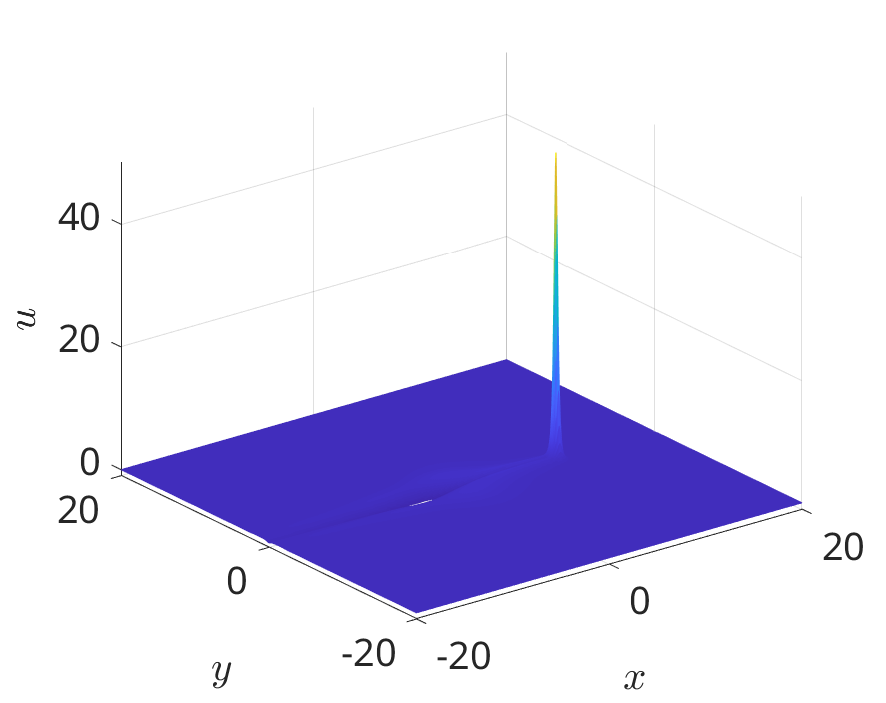}
\caption{Evolution of the double perturbed solitons from \eqref{E:doubleQ} at $t=1$ on the left, $t=2.5$ in the middle, $t= 3.5$ on the right.}
\label{F:doubleQ-1}
\end{figure}

We examine the situation similar to the one in the previous subsection, replacing the Gaussians with solitons. More precisely, we take two-bump initial condition 
\begin{equation}\label{E:doubleQ}
u = 0.95Q(x, y-\pi/2,z) + 0.95Q(x, y+\pi/2,z).
\end{equation}
These are two overlapping perturbed solitons, with the peaks separated by $\pi$ in the $y$-direction. 
Each component has mass $0.95^2M[Q]=0.9025M[Q]$, and the positive overlap makes the total mass larger than $1.8 \, M[Q]$. 
Numerically computing the energy of this double-soliton, we get 
$$
E[u_0] \approx -12.5910,
$$
which is negative (unlike the double-Gaussian case in the previous example). We note that each individual bump has positive energy
$$
E[0.95Q] = \tfrac12 \, 0.95^2 \|\nabla Q\|^2_{L^2} - \tfrac{10}3 0.95^{3/{10}} \|Q\|^{10/3}_{L^{10/3}} = \tfrac34 (0.95^2 - 0.95^{10/3}) M(Q) \approx 2.85,
$$ 
where we used \eqref{E:Grad-Pot} and the value of the mass from \eqref{E:L2}. If the interaction would be insignificant between two bumps, then the energy would remain positive (the energy of two separate soliton bumps would be around 5.7). However, the energy is quite negative, so the overlap removes a significant chunk of energy, which in a sense influences the global behavior, which we show in Figure \ref{F:doubleQ-1}.
This is due to the slower decay of the ground state (as $\frac{e^{-r}}{r}$) compared to the Gaussian (as $e^{-r^2}$), thus, influencing interaction of bumps more on the same separation interval.   

\smallskip

\begin{figure}[!htb]
\includegraphics[width=0.32\hsize]{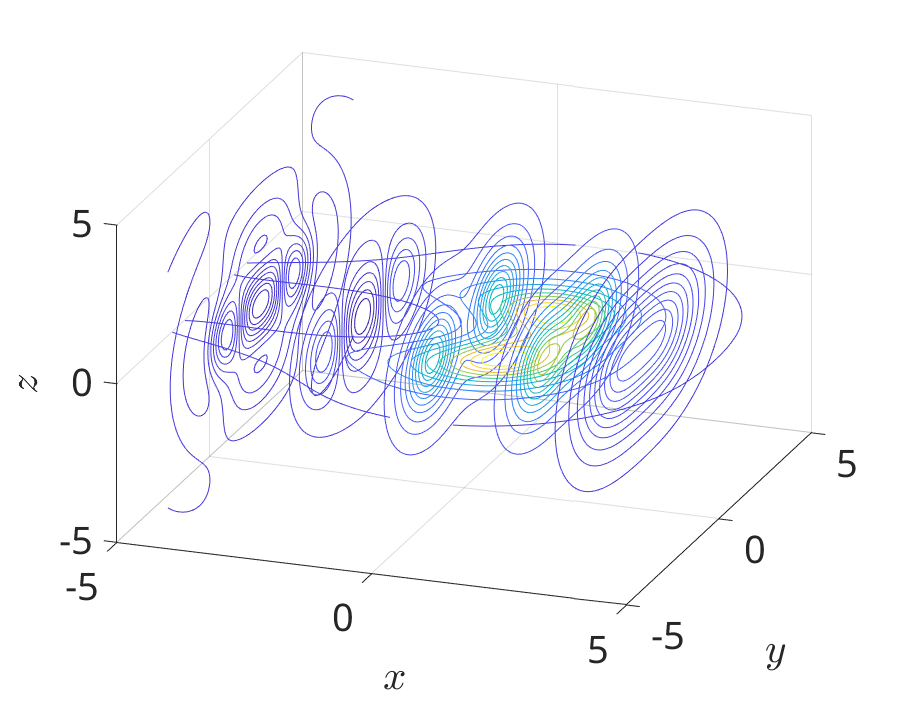}
 \includegraphics[width=0.32\hsize]{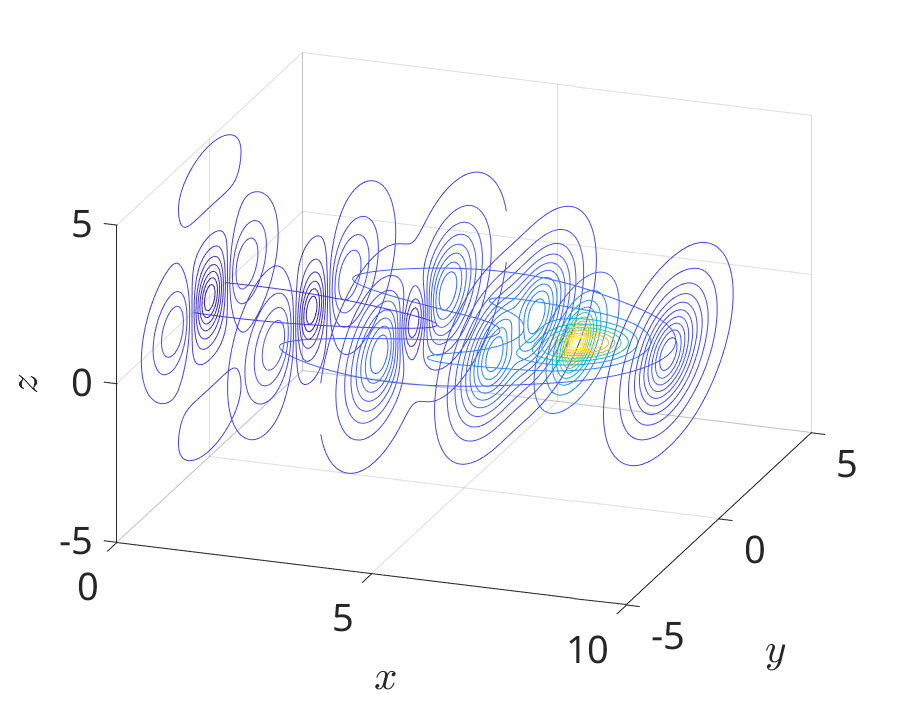}
  \includegraphics[width=0.32\hsize]{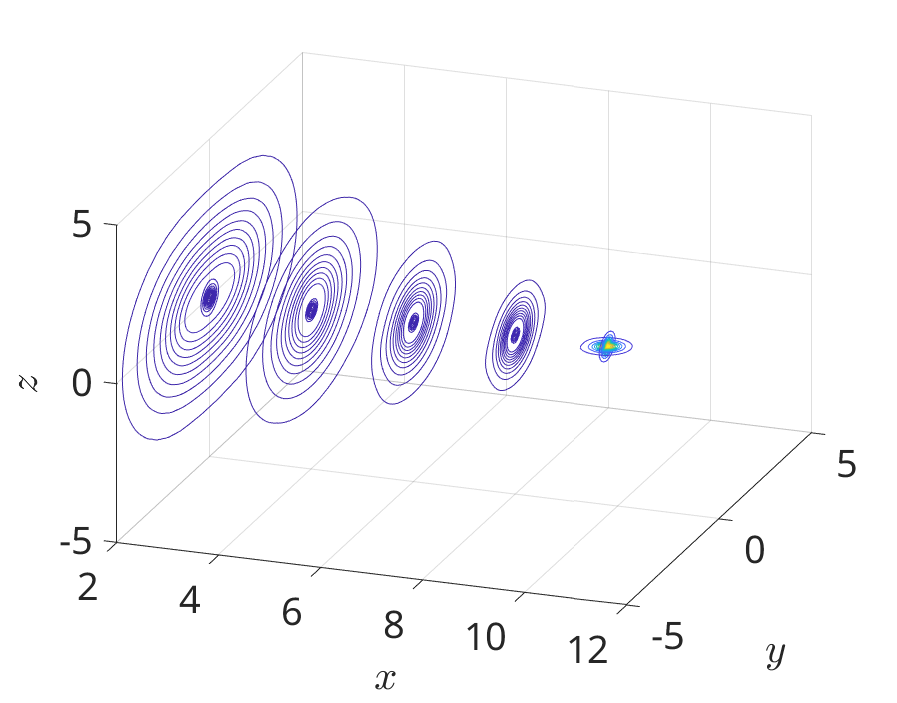}
\caption{Contour plots for the evolution of double perturbed soliton 
from \eqref{E:doubleQ} at $t=1$ on the left, $t=2.5$ in the middle, 
$t= 3.5$ on the right. To improve visualisation we zoom into the 
region near  the mass distribution and present a contour plot on 
planes, for even values of $x$ plus one through the maximum of the 
solution. There are 10 isolines on each plane. Note that the direction of the $x$-axis is reversed compared to the standard to improve presentation. We use the periodicity in $x$ as before.} 
\label{F:doubleQ-2}
\end{figure}

\begin{figure}[!htb]
\includegraphics[width=0.32\hsize]{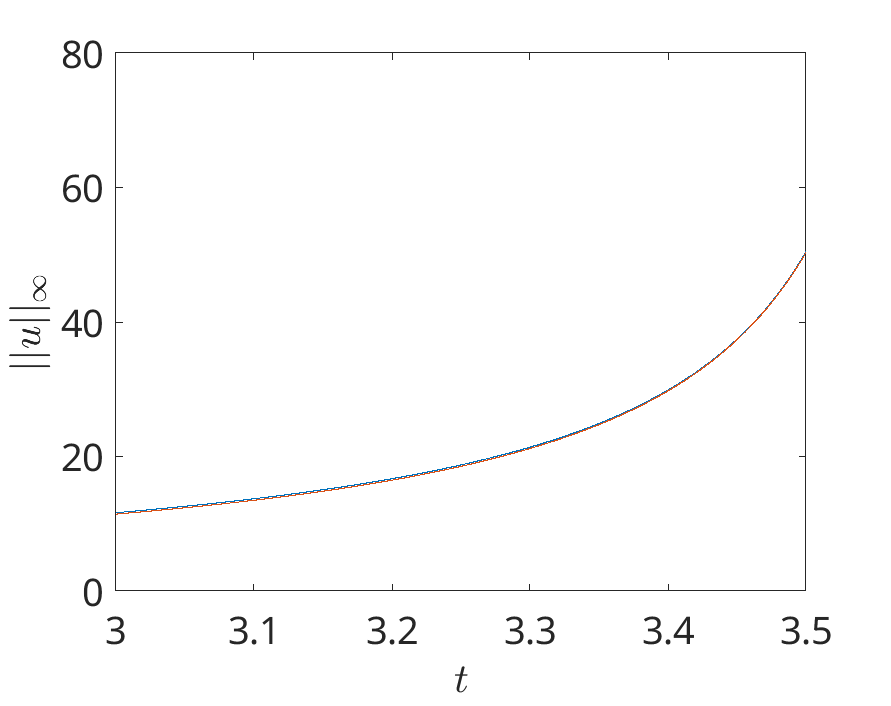}
\includegraphics[width=0.32\hsize]{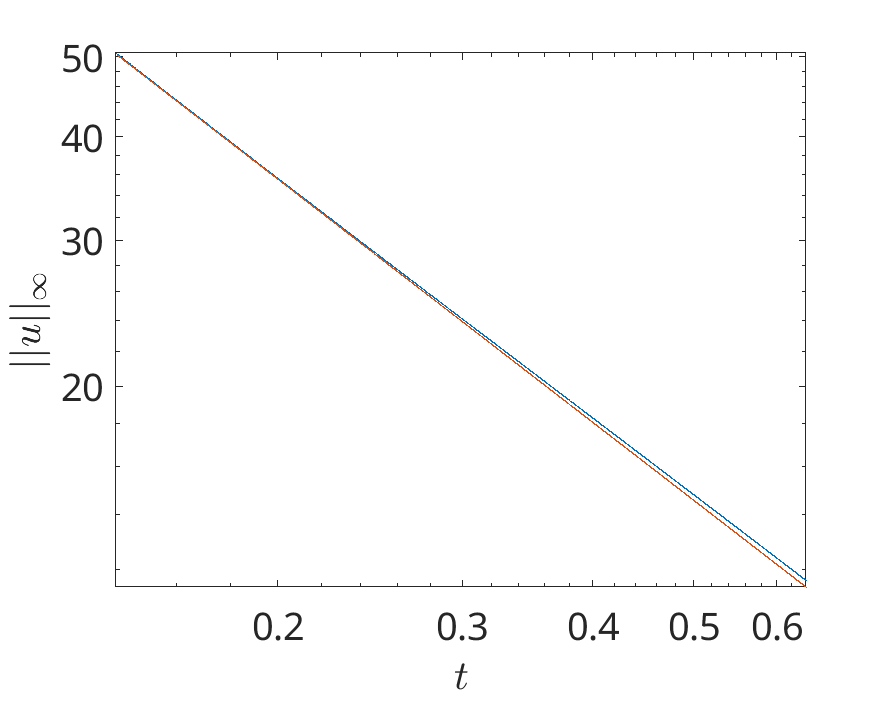}
\includegraphics[width=0.32\hsize]{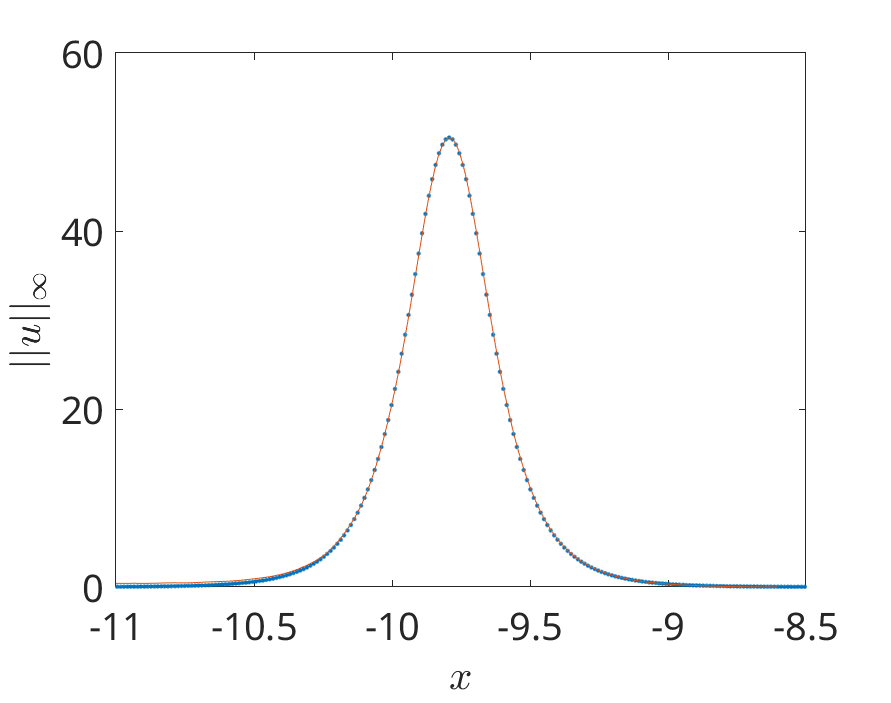}
\caption{On the left a section of the profile of the solution along $x$ in red and a rescaled soliton in blue. In the middle the $\|u\|_{L^\infty}$ norm in blue and a blowup rate following $b(t_c-t)^a$ in red. On the right the same two curves in a log-log plot. We have found $t_c = 3.64$, $a = -0.975$, $b = 7.4$.}
\label{F:doubleQ-3}
\end{figure}

Our simulations show that with time the two perturbed soliton bumps merge and form a single peak, which then continues to concentrate indicating the blow-up in finite time, see the time evolution in Figure~ \ref{F:doubleQ-1}. Their contour plots are shown in Figure ~\ref{F:doubleQ-2} to further indicate the formation of radiation and the concentration as singularity forms.

We track the time dependence of the $L^\infty$ norm in Figure~\ref{F:doubleQ-3} and estimate the parameters of the power law blow-up, noting that the blow-up rate here matches closer the predicted rate in \S \ref{S:3Drate}: here the fitted exponent $a=-0.975$ is close to the formal prediction for the $L^\infty$ norm $-3/2 \, \beta_{3d} \approx -0.967$.  We also show the profile of the solution at the final computational time $t=3.5$ (with the estimated blow-up time $t_c=3.64$) matching well with a rescaled ground state, see right plot in Figure ~\ref{F:doubleQ-3}. 


\section{Conclusions}\label{S:conclusions}

In this paper we study the three-dimensional $L^2$-critical ZK equation with the fractional nonlinearity $|u|^{4/3}u$ ($\equiv u^{7/3}$ for real-valued solutions $u$) from formal analytical and full 3D computational points of view.  Extending the two-dimensional approach used for the critical ZK equation in 2D in \cite{Gong} we obtain formulas for the blow-up correction coefficient $\theta_{3d} \approx 1.448$, and the resulting formal prediction for the blow-up rate $\beta_{3d} \approx 0.644$, where
$$
\lambda(t)\sim C(T-t)^{\beta_{3d}}, \qquad  \|\nabla u(t)\|_{L^2}\sim C(T-t)^{-\beta_{3d}}, \qquad  \|u(t)\|_{L^\infty}\sim C(T-t)^{-3/2\, \beta_{3d}}. 
$$
We also provide an alternative way to compute these coefficients using the fact that the ground state is radial, thus, an application of the projection-slice theorem and radial Fourier transform reduce the computation of $\theta$'s to a one-dimensional sine-transform quotient, which could be  inexpensive to evaluate once the radial ground state is known. 
We provide a formal analysis for blow-up rate in any dimension for the $L^2$-critical ZK, and give values for the 2D and 4D cases.  

Computationally, we implement the full three-dimensional GPU simulations, which show several robust qualitative features. Studying perturbation of single-bump initial data such as the ground state and the negative-energy Gaussian, we show that the solutions develop a sharply concentrating core whose profile is well approximated by a rescaled ground state on the computed time interval. The simulations also show a radiative part propagating in a cone region opposite to the solitary-wave direction.  The double-soliton experiment produces a merged concentrating core and an $L^\infty$ exponent close to the formal prediction.  By contrast, the double-Gaussian data have total mass above $M[Q]$ but disperses part of it into radiation and then stabilize, on the computed time interval. This comparison shows that total mass alone does not determine the observed dynamics and that localization, profile shape, overlap, and energy must be taken into account.

The rate computations show some limitations. The fitted exponents for the $1.1Q$ and single-Gaussian runs provide evidence for the blow-up rates in the pre-asymptotic regime, while showing the leading profile being a ground state. This can be explained due to the final computational time that we can reach at the moment. In the two-bump blow-up dynamics, our simulations show a closer value to the predicted blow-up rate again with the leading blow-up profile being a rescaled ground state.   Reaching a definitive rate comparison will require significantly higher computational power as well as significantly larger memory together with systematic refinements, adaptive meshes, and a fitting protocol that allows to quantify the singularity parameters more precisely.



\end{document}